\documentclass[a4paper,10pt,twoside]{article}
\usepackage{pst-fill,pst-grad}
\usepackage{textcomp}
\usepackage[english]{babel}
\usepackage[utf8x]{inputenc}
\usepackage{graphicx}
\usepackage{amsmath}
\usepackage{float}
\usepackage{fancyhdr}
\usepackage[matrix,arrow,curve]{xy}
\usepackage{pstricks} 
\usepackage{amsmath,amsfonts,verbatim,afterpage,theorem,euscript,mathrsfs,amssymb}
\usepackage{amsfonts}
\usepackage{amssymb}
\usepackage{array}
\usepackage{dsfont}
\usepackage{hyperref}
\newtheorem{Definition}{Definition}[section]
\newtheorem{Proposition}{Proposition}[section]
\newtheorem{Lemme}{Lemma}[section]
\newtheorem{Theoreme}{Theorem}

\newtheorem{Remarque}{Remark}[section]
\title{\bf A molecular method applied to a non-local PDE in stratified Lie groups}
\author{Diego Chamorro\footnote{Laboratoire d'Analyse et de Probabilit\'es, Universit\'e d'Evry Val d'Essonne, 23 Boulevard de France, 91037 Evry Cedex - France, diego.chamorro@univ-evry.fr}}
\begin{document}
\maketitle
\begin{scriptsize}
\abstract{In this article we study a transport-diffusion equation in the framework of the stratified Lie groups. For this equation we will study the existence of the solutions, a maximum principle, a positivity principle and H\"older regularity. }\\[3mm]
\end{scriptsize}
\section{Introduction}
The Euclidean space $\mathbb{R}^n$ has many useful properties and all the standard tools are adapted in a very natural way to its structure. However, if we slightly change the setting by considering small modifications in the general structure, some of the usual and standard technics used in the analysis of the PDEs can possibly lose all their interest as the computations becomes potentially more involved. \\

In this article, we will consider a different framework and we will replace the space $\mathbb{R}^n$ by stratified Lie groups. These groups are a generalization of $\mathbb{R}^n$ with a different dilation structure and a different group law. For example, if $x=(x_{1}, x_{2}, x_{3})$ is an element of $\mathbb{R}^{3}$, we can fix a dilation by writing 
$\delta_{\alpha}[x]=(\alpha x_{1 }, \alpha x_{2}, \alpha^{2 } x_{3})$ for $\alpha>0$. Then, the adapted group law with respect to this dilation (in the sense that $\delta_{\alpha}[x\cdot y] = \delta_{\alpha}[x]\cdot \delta_{\alpha}[y] $) is given by  
\begin{equation*}
 x\cdot y=(x_{1}, x_{2}, x_{3})\cdot(y_{1}, y_{2}, y_{3})=(x_{1}+y_{1}, x_{2}+y_{2}, x_{3}+y_{3}+\frac{1}{2}(x_{1}y_{2}-y_{1}x_{2 })). 
\end{equation*} 
This triplet $(\mathbb{R}^{3}, \cdot,\delta )$ corresponds to the Heisenberg group $\mathbb{H}^{1}$ which is the first non-trivial example of a stratified Lie group. Note in particular that $x \cdot y \neq y \cdot x$ so this group is no longer abelian. Remark also that although the topological dimension of $\mathbb{H}^{1}$ is $n=3$, the homogeneous dimension with respect to the dilation $\delta_\alpha$ given above is $N=4$ and we will see how this particular fact gives a special flavor to the computations. Furthermore it can be shown that the structure of $\mathbb{H}^{1}$ is completely different from the usual one of $\mathbb{R}^{3}$, thus these seemingly small changes in the inner structure leads to deep modifications that must be taken into account, see \cite{Stein2} for more details.\\ 

Once such a dilation and group structure is fixed, all the corresponding objects (such as operators and functional spaces) are defined in order to follow it and this imposes some specific problems that are not merely technical. Take for example the fact that if the underlying group structure is not commutative, we can not assure in all generality the helpful identity $f\ast g = g\ast f $ for the convolution of two functions. Another example is raised with the use of the Fourier transform: even if we can explicitely use it in some cases to perform specific calculations, these can occasionally be quite complex.\\

The aim of this article is to adapt to the setting of stratified Lie groups a recent method developped in \cite{Chame} and to show its robustness for the study of a transport-diffusion equation. We will see then how this method, based on the H\"older-Hardy spaces duality and on the characterization of the Hardy spaces in terms of molecules, can be successfully generalized into this family of Lie groups.\\ 

We are going to work with the partial diferential equation of the following form:
\begin{equation}\label{Equation0}
\begin{cases}
\partial_t \theta-\nabla\cdot(v\,\theta)+\mathcal{J}^{1/2}\theta=0,\\[5mm]
\theta(x,0)= \theta_0(x),\\[5mm]
\mbox{with }\; div(v)= 0\; \mbox{ and } t\in [0, T].
\end{cases}
\end{equation}
Here $\theta$ is a function over a stratified Lie group $\mathbb{G}$, $v$ is a divergence-free velocity field that belongs to the space of functions of local bounded mean oscillations $bmo(\mathbb{G})$ and $\mathcal{J}^{1/2}$ stands for the square-root of a sub-Laplacian of the group $\mathbb{G}$. We will give all the definitions in section \ref{Definition} below.\\

This type of transport-diffusion equations is a generalization of a well-known equation from fluid dynamics: in the Euclidean setting $ \mathbb{R}^2$ if $\mathcal{J}^{1/2}=(-\Delta)^{1/2}$ is the fractional Laplacian and if $v=(-R_2\theta, R_1\theta)$ where $R_{1,2}$ are the Riesz Transforms defined in the Fourier level by $ \widehat{R_j\theta}(\xi)=-\frac{i\xi_j}{|\xi|}\widehat{\theta}(\xi)$ for $j=1,2$, we obtain the quasi-geostrophic equation $(QG)_{1/2}$. See \cite{Caffarelli}, \cite{Cordoba}, \cite{CW1}, \cite{CW}, \cite{Marchand1} and the references there in for more details.\\

It is worth noting in this type of equations that there is a competition between the transport term and the diffusion one which is given by the fractional power of the Laplacian. In the setting studied with equation (\ref{Equation0}), each one of these two terms have formally the same effect and with the useful hypotesis $div(v)=0$ we will prove that it is possible to obtain a small gain of regularity. 
\subsection*{Presentation of the results}

Although many tools are available in the setting of stratified Lie groups, these objects are mainly studied for their own sake and, to the best of our knowledge, the treatment of equation (\ref{Equation0}) in this framework is new. For this reason we will prove step by step
some useful properties of the solutions to this equation with the next theorems.\\ 

In what follows, we will always assume that $\mathbb{G}$ is a stratified Lie group of homogeneous dimension $N\geq 3$.\\

\begin{Theoreme}[Existence and uniqueness for $L^p$ initial data]\label{Theo0}
If $\theta_0\in L^{p}(\mathbb{G})$ with $2\leq p\leq +\infty$ is an initial data, then equation (\ref{Equation0}) has a unique weak solution $\theta\in L^\infty([0,T]; L^p(\mathbb{G}))$.
\end{Theoreme}

\begin{Theoreme}[Maximum and Positivity Principle]\label{Theo1}
\begin{itemize}
\item[]
\item[1)] If $\theta$ is a smooth solution of equation (\ref{Equation0}), then we have for $1\leq p \leq +\infty$ the following \emph{maximum principle}: 
$$\|\theta(\cdot, t)\|_{L^p}\leq \|\theta_0\|_{L^p}.$$
\item[2)] Let $N<p \leq +\infty$ and $M>0$ a constant, if the initial data $\theta_0\in L^p(\mathbb{G})$ is such that $0\leq \theta_0\leq M$ then the weak solution of equation (\ref{Equation0}) satisfies the \emph{positivity principle}: we have $0\leq \theta(x,t)\leq M$ for all $t\in [0, T]$. 
\end{itemize}
\end{Theoreme}
Our main theorem is the following one:
\begin{Theoreme}[Hölder regularity]\label{Theo3}
Fix a small time $T_0>0$ and consider $\theta_0\in L^{p}(\mathbb{G})$ with $2\leq p\leq +\infty$. If $\theta(x,t)$ is a solution for the equation (\ref{Equation0}) with associated initial data  $\theta_0$, then for all time $T_0<t<T$, we have that $\theta(\cdot,t)$ belongs to the H\"older space $\mathcal{C}^{\gamma}(\mathbb{G})$ with $0<\gamma<1$. 
\end{Theoreme}

Some remarks are in order. Theorems \ref{Theo0} and \ref{Theo1} will be proven with rather classical ideas using the appropriate tools: the important point here is the application of well suited characterization of functional spaces in the framework of stratified Lie groups.\\

The lower bound $2\leq p$ for the existence of the solutions and for the maximum principle is mainly due to technicalities and it should be possible to treat the case $1\leq p <2$. The condition $N< p$ for the positivity principle is also technical. However, the cases $2\leq p$ and $N< p$ are enough for our purposes. \\

For proving Theorem \ref{Theo3} we will adapt a duality based technic recently developped in $\mathbb{R}^n$ to the setting of stratified Lie groups. This method relies in the use of molecular Hardy spaces which have been broadly studied in many different situations: in particular for stratified groups we have at our diposal an explicit characterization with \textit{molecules} (see \cite{DeMichele}).\\

The plan of the article is the following: in section \ref{Definition} we give the principal definitions, in the section \ref{SecExiUnic} we study existence and uniqueness of solutions with initial data in $L^p$ with $2\leq p<+\infty$. In this section we will also prove the maximum principle. Section \ref{Sect_PrincipeMax} is devoted to a positivity principle that will be useful in our proofs and section \ref{SecLinfty} studies existence of solution with $\theta_0\in L^\infty$. In section \ref{SeccHolderRegularity} we study the H\"older regularity of the solutions of equation (\ref{Equation0}) by a duality method.
\section{Definitions}\label{Definition}
In this section we recall some basic facts about stratified Lie groups, for further information see \cite{Folland2}, \cite{Varopoulos},\cite{Stein2}  and the references given there. \\

A \textit{homogeneous group} $\mathbb{G}$ is the data of $\mathbb{R}^{n}$ equipped with a structure of Lie group and with a family of dilations 
which are group automorphisms and it will be denoted by $\mathbb{G}=(\mathbb{R}^{n}, \cdot, \delta)$. We will always suppose that the origin is the identity. We define the \textit{dilations} by fixing integers $(a_{i})_{1\leq i\leq n}$ such that $1=a_{1}\leq... \leq a_{n}$ and by writing:  
\begin{eqnarray*}
\delta_{\alpha }:  \mathbb{R}^{n } & \longrightarrow & \mathbb{R}^{n } \\
x &\longmapsto & \delta_{\alpha}[x]=(\alpha^{a_{1}}x_{1},...,\alpha^{a_{n}}x_{n})\nonumber \qquad (\alpha>0).
\end{eqnarray*} 
The \emph{homogeneous dimension} with respect to this dilation is given by $N=\displaystyle{\sum_{1\leq i\leq n}}a_{i }$. Observe that $N\geq n$ since $a_{i}\geq 1$ for all $i=1,...,n$. We will say that a function on $\mathbb{G}\setminus \{0\}$ is \textit{homogeneous of degree} $\lambda \in \mathbb{R}$ if $f(\delta_{\alpha}[x])=\alpha^{\lambda}f(x)$ for all $\alpha>0$.  In the same way, we will say that a differential operator $D$ is homogeneous of degree $\lambda$ if $D(f(\delta_{\alpha}[x]))=\alpha^{\lambda}(Df)(\delta_{\alpha}[x])$ for all $f$ in the operator's domain.  
In particular, if $f$ is homogeneous of degree $\lambda$  and if $D$ is a differential operator of degree $\mu$, then $Df$ is homogeneous of degree $\lambda-\mu$. \\ 

In homogeneous groups the Lebesgue measure $dx$ is bi-invariant and coincides with the Haar measure. For any subset $E$ of $\mathbb{G}$ we note its measure as $|E|$ and we define the $L^p$ norms in the usual way:
\begin{equation*}
\|f\|_{L^p}=\left(\int_{\mathbb{G}}|f(x)|^pdx\right)^{1/p}\quad  \mbox{ if } 1\leq p <+\infty \quad  \mbox{ and } \quad  \|f\|_{L^\infty}=\underset{x\in \mathbb{G}}{\sup \, ess} |f(x)|\quad \mbox{ if } p=+\infty.
\end{equation*}
The convolution of two functions $f$ and $g$ on $\mathbb{G}$ is then defined by
\begin{equation*} 
f\ast g(x)=\int_{\mathbb{G}}f(y)g(y^{-1}\cdot x)dy=\int_{\mathbb{G}}f(x\cdot y^{-1})g(y)dy,  \quad x\in \mathbb{G}.
\end{equation*}
We also have the useful Young's inequalities:
\begin{Lemme}
If $1\leq p, q, r\leq +\infty$ such that $1+\frac{1}{r}=\frac{1}{p}+\frac{1}{q}$.  
If $f\in L^{p}(\mathbb{G})$ and $g\in L^{q}(\mathbb{G})$, then $f\ast g \in L^{r}(\mathbb{G})$ and 
\begin{equation*}
\|f\ast g\|_{L^r}\leq \|f\|_{L^p}\|g\|_{L^q}. 
\end{equation*} 
\end{Lemme} 
A proof is given in \cite{Folland2}.\\

For a homogeneous group $\mathbb{G}=(\mathbb{R}^{n}, \cdot, \delta)$ we consider now its Lie algebra $\mathfrak{g}$ whose elements can be conceived in two different ways: as \textit{left}-invariant vector fields or as \textit{right}-invariant vector fields. The left-invariant vectors fields $(X_j)_{1\leq j\leq n}$ are determined by the formula
\begin{equation*}
(X_{j}f)(x)=\left.\frac{\partial f(x\cdot y)}{\partial y_{j}}\right|_{y=0}=\frac{\partial f}{\partial x_{j}}+\sum_{j<k}q^{k}_{j}(x)\frac{\partial f}{\partial x_{k}} 
\end{equation*} 
where $q^{k}_{j}(x)$ is a homogeneous polynomial of degree $a_{k}-a_{j}$ and $f$ is a smooth function on $\mathbb{G}$. By this formula one deduces easily that these vectors fields are homogeneous of degree $a_{j}$: $X_{j}\left(f(\alpha x)\right)=\alpha^{a_{j}}(X_{j}f)(\alpha x)$. We will note $(Y_{j})_{1\leq j\leq n}$ the right invariant vector fields defined in a totally similar way:
$$(Y_{j}f)(x)=\left.\frac{\partial f(y\cdot x)}{\partial y_{j}}\right|_{y=0}$$
A homogeneous group $\mathbb{G}$ is  \emph{stratified} if its Lie algebra $\mathfrak{g}$ breaks up into a sum of linear subspaces  
$\mathfrak{g}=\displaystyle{\bigoplus_{1\leq j\leq k }} E_{j}$ such that $E_{1}$ generates the algebra $\mathfrak{g}$ and $[E_{1}, E_{j}]=E_{j+1}$ for $1\leq j < k$ and $[E_{1}, E_{k}]=\{0\}$ and $E_{k}\neq\{0\}$, but $E_{j}=\{0\}$ if $j>k$. Here $[E_{1}, E_{j}]$ indicates the subspace of $\mathfrak{g}$ generated by the elements $[U, V]=UV-VU$ with $U\in E_{1}$ and $V\in E_{j}$. The integer $k$ is called the \emph{degree} of stratification of $\mathfrak{g}$. For example, on Heisenberg group $\mathbb{H}^1$, we have $k=2$ while in the Euclidean case $k=1$.  \\

We will suppose henceforth that $\mathbb{G}$ is \textbf{stratified}. Within this framework, if we fix the vectors fields $X_{1},...,X_{m}$ such that $a_{1}=a_{2}=\ldots=a_{m}=1$ $(m<n)$, then the family $(X_{j})_{1\leq j\leq m}$ is a base of $E_{1}$ and generates the Lie algebra of $\mathfrak{g}$, which is precisely the H\"ormander's condition (see \cite{Folland2} and \cite{Varopoulos}).\\

To the family $(X_{j})_{1\leq j\leq m}$ is associated the Carnot-Carath\'eodory distance $d$ which is left-invariant and compatible with the topology on $\mathbb{G}$ (see \cite{Varopoulos} for more details). For any $x\in \mathbb{G}$ we will note $\|x\|=d(x,e)$ and $d(x,y)=\|x\cdot y^{-1}\|$. For $r>0$ we form the balls by writing $B(x,r)=\{y\in \mathbb{G}: d(x,y)<r\}$.\\

For any multi-index $I=(i_{1},...,i_{n})\in \mathbb{N}^{n}$, one defines $X^{I}$ by $X^{I}=X_{1}^{i_{1}}\dots X_{n}^{i_{n}}$ and $Y^{I}$ by $Y^{I}=Y_{1}^{i_{1}}\dots Y_{n}^{i_{n}}$. We note $|I|= i_{1}+\ldots+i_{n}$ the order of the derivation $X^I$ or $Y^I$ and $d(I)=a_{1}i_{1}+\ldots+a_{n}i_{n}$ the homogeneous degree of this one. \\  

For $\varphi, \psi \in \mathcal{C}^{\infty}_{0}(\mathbb{G})$ we have the equality 
\begin{equation*} 
\int_{\mathbb{G}}\varphi(x)(X^{I}\psi)(x)dx=(-1)^{|I|}\int_{\mathbb{G}}(X^{I}\varphi)(x)\psi(x)dx.  
\end{equation*} 
The interaction of operators $X^{I}$ and $Y^{I}$ with convolutions is clarified by the following identities:  
\begin{equation*}
X^{I}(f*g)=f*(X^{I}g), \qquad Y^{I}(f*g)=(Y^{I}f)*g, \qquad (X^{I}f)*g=f*(Y^{I}g). 
\end{equation*} 
We will say that a function $f \in \mathcal{C}^{\infty}(\mathbb{G})$ belongs to the Schwartz class $\mathcal{S}(\mathbb{G})$ if the following semi-norms are bounded for all $k\in \mathbb{N}$ and any multi-index $I$:  $N_{k,I}(f)=\underset{x\in \mathbb{G}}{\sup } \, (1+|x|)^{k }|X^{I}f(x)|$. 
\begin{Remarque}
To characterize the Schwartz class $\mathcal{S}(\mathbb{G})$ we can replace vector fields $X^I$ in the semi-norms $N_{k,I}$ above by right-invariant vector fields $Y^I$. For a proof of these facts and for further details see \cite{Folland2}.
\end{Remarque}
We define now the \textit{gradient} on $\mathbb{G}$ from vectors fields of homogeneity degree equal to one by fixing $\nabla = (X_{1},...,X_{m })$.  This operator is of course left invariant and homogeneous of degree $1$. The length of the gradient is given by the formula $|\nabla f|= \left((X_{1}f)^{2}+... +(X_{m}f)^{2 } \right)^{1/2}$. Once the gradient is fixed, we will work with the following sub-Laplacian:
\begin{equation*}
\mathcal{J}=\nabla^{*}\nabla=-\sum_{j=1}^{m}X^{2}_{j} 
\end{equation*} 
which is a positive self-adjoint, hypo-elliptic operator since $(X_j)_{1\leq j\leq m}$ satisfies the H\"ormander's condition. Its associated \textit{heat operator} on $\mathbb{G}\times]0, +\infty[$ is given by $\partial_{t}+\mathcal{J}$. We recall now some well-known properties:
\begin{Theoreme} The semi-group $H_{t}=e^{-t\mathcal{J}}$  admits a convolution kernel $H_{t}f=f\ast h_{t}$ where $h_{t}(x)=h(x,t) \in \mathcal{C}^{\infty}(\mathbb{G}\times]0, +\infty[)$ is the heat kernel which satisfies the following points:  
\begin{itemize}
\item $h(x, t)=h(x^{-1}, t)$, $h(x, t)\geq 0$ and $\displaystyle{\int_{\mathbb{G}}}h(x, t)dx=1$, 
\item $h_{t}$ has the semi-group property:  $h_{t}\ast h_{s}=h_{t+s}$ for $t, s>0$, 
\item $h(\delta_{\alpha}[x], \alpha^{2}t)=\alpha^{-N}h(x, t)$, 
\item For every $t>0$, $x\mapsto h(x, t)$ belongs to the Schwartz class in $\mathbb{G}$. 
\end{itemize} 
\end{Theoreme} 
For a detailed proof of these and other important facts concerning the heat semi-group see \cite{Folland2} and \cite{Saka}.\\

In order to define regularity measuring spaces, we will need to deal with fractional powers of the sub-Laplacian. For $s>0$ we define $\mathcal{J}^{s/2}$ using the spectral resolution of the sub-Laplacian: $\mathcal{J}^{s/2}=\displaystyle{\int_{0}^{+\infty}}\lambda^{s/2}dE(\lambda)$, and we have that the operator $\mathcal{J}^{s/2}$ is homogeneous of degree $s$.\\

In this article we will mainly work with the square root of the sub-Laplacian $\mathcal{J}^{1/2}$ and we will need the following equivalent characterization:
$$\mathcal{J}^{1/2} f(x)\simeq v.p. \int_{\mathbb{G}}\frac{f(x)-f(y)}{\|x\cdot y^{-1}\|^{N+1}}dy$$
This definition is related to the properties of the Poisson semi-group and its kernel. See more details in \cite{Folland} and \cite{Folland1}.\\

We define now the Sobolev space $W^{s,p}(\mathbb{G})$ with $1<p<+\infty$ and $s>0$ by the norm
\begin{equation*}
\|f\|_{W^{s,p}}=\|f\|_{L^p}+\|\mathcal{J}^{s/2}f\|_{L^p}
\end{equation*}
It's homogeneous version is given by $\|f\|_{\dot W^{s,p}}=\|\mathcal{J}^{s/2}f\|_{L^p}$. For Besov spaces  $B^{s,p}_q(\mathbb{G})$ we have, for $0<s<1$ and  $1\leq p,q\leq +\infty$:
\begin{equation*}
\|f\|_{B^{s,p}_q}= \|f\|_{L^p}+\|f\|_{\dot B^{s,p}_q}
\end{equation*}
where we noted 
\begin{equation*}
\|f\|_{\dot B^{s,p}_q}=\left(\int_{\mathbb{G}} \left[\int_{\mathbb{G}} \frac{|f(x\cdot y)-f(x)|^p}{\|y\|^{(N+sq)\frac{p}{q}}}dx\right]^{q/p}dy\right)^{1/q}.
\end{equation*}
To finish, we define the space $bmo(\mathbb{G})$ as the space of locally integrable functions such that 
$$\underset{|B|\leq 1}{\sup}\frac{1}{|B|}\int_{B}|f(x)-f_B|dx<M \qquad\mbox{ and }\qquad \underset{|B|>1}{\sup}\frac{1}{|B|}\int_{B}|f(x)|dx<M\qquad \mbox{for a constant } M;$$
 where we noted $B(R)$ a ball of radius $R>0$ and $f_B=\frac{1}{|B|}\displaystyle{\int_{B(R)}}f(x)dx$. The norm $\|\cdot\|_{bmo}$ is then fixed as the smallest constant $M$ satisfying these two conditions. 

\section{Existence and uniqueness with $L^p$ initial data and Maximum Principle.}\label{SecExiUnic}
In this section we will study existence and uniqueness for weak solution of equation (\ref{Equation0}) with initial data $\theta_0\in L^p(\mathbb{G})$ where $p\geq 2$. We will start by considering an approximation of this equation and we will prove existence and uniqueness for this system. To pass to the limit we will need a further step that is a consequence of the maximum principle.\\

We begin our study with the following approximation of equation (\ref{Equation0}):
\begin{equation}\label{SistApprox}
\left\lbrace
\begin{array}{l}
\partial_t \theta+ \nabla\cdot(v_{\varepsilon}\;\theta)+\mathcal{J}^{1/2}\theta= - \varepsilon  \mathcal{J}\theta \qquad \qquad (\varepsilon>0)\\[5mm]
\theta(x,0)=\theta_0(x)\\[5mm]
div(v)=0 \quad \mbox{ and } v\in L^{\infty}([0,T];  L^{\infty}(\mathbb{G})).
\end{array}
\right.
\end{equation}
where $v_{\varepsilon}$ is defined by $v_{\varepsilon}=v\ast \omega_{\varepsilon}$ with $\omega_{\varepsilon}(x)=\varepsilon^{-N}\omega(\delta_{\varepsilon^{-1}}[x])$ and $\omega\in \mathcal{C}^{\infty}_0(\mathbb{G})$ is a function such that $\displaystyle{\int_{\mathbb{G}}}\omega(x)dx=1$. 
\begin{Remarque}
It is equivalent to consider $-v$ instead of $v$, thus for simplicity we fix velocity's sign as in equation (\ref{SistApprox}) above. The same proofs are valid for equation (\ref{Equation0}). 
\end{Remarque}
Observe that we fixed here the velocity $v$ such that $v\in L^{\infty}([0,T'];  L^{\infty}(\mathbb{G}))$. This is not very restrictive since we have the following lemma:
\begin{Lemme}\label{TheoApproxbmo}
Let $f$ be a function in $bmo(\mathbb{G})$. For $k\in \mathbb{N}$, define $f_k$ by
\begin{equation*}
f_k(x)=\left\lbrace
\begin{array}{rll}
-k & \mbox{if} & f(x)\leq -k \\[2mm]
f(x) & \mbox{if} & -k\leq f(x)\leq k \\[2mm]
k & \mbox{if} & k\leq f(x).
\end{array}
\right.
\end{equation*}
Then $(f_k)_{k\in \mathbb{N}}$ converges $\ast$-weakly to $f$ in $bmo(\mathbb{G})$.
\end{Lemme}
A proof of this lemma can be found in \cite{Stein2}.\\

Note now that the problem (\ref{SistApprox}) admits the following equivalent integral representation:
\begin{equation}\label{FormIntegr}
\theta(x,t)=e^{-\varepsilon t \mathcal{J}}\theta_0(x)-\int_{0}^t e^{-\varepsilon (t-s) \mathcal{J}}\nabla \cdot(v_\varepsilon\; \theta)(x,s)ds-\int_{0}^t e^{-\varepsilon (t-s) \mathcal{J}}\mathcal{J}^{1/2} \theta(x,s)ds,
\end{equation}
In order to prove Theorem \ref{Theo0}, we will first investigate a local result with the following theorem where we will apply the Banach contraction scheme in the space $L^{\infty}([0,T]; L^{p}(\mathbb{G}))$ with the norm $\|f\|_{L^\infty (L^p)}=\displaystyle{\underset{t\in [0,T]}{\sup}}\|f(\cdot, t)\|_{L^p}$.
\begin{Theoreme}[Local existence for viscosity solutions]\label{TheoPointFixe}
Let $1\leq p<+\infty$ and let $\theta_0$ and $v$ be two functions such that $\theta_0\in L^p(\mathbb{G})$, $div(v)=0$ and $v\in L^{\infty}([0,T']; L^{\infty}(\mathbb{G}))$. If the initial data satisfies $\|\theta_0\|_{L^p}\leq K$ and if $T'$ is a time small enough, then (\ref{FormIntegr}) has a unique solution $\theta \in L^{\infty}([0,T']; L^{p}(\mathbb{G}))$ on the closed ball $\overline{B}(0,2K)\subset L^{\infty}([0,T']; L^{p}(\mathbb{G}))$. 
\end{Theoreme}
\textit{\textbf{Proof of Theorem \ref{TheoPointFixe}.}} We construct a sequence of functions in the following way
$$\theta_{n+1}(x,t)=e^{-\varepsilon t \mathcal{J}}\theta_0(x)-\int_{0}^t e^{-\varepsilon (t-s) \mathcal{J}}\nabla \cdot(v_\varepsilon\; \theta_n)(x,s)ds-\int_{0}^t e^{-\varepsilon (t-s) \mathcal{J}}\mathcal{J}^{1/2}\theta_n(x,s)ds,$$ 
and we take the $L^\infty L^p$-norm of this expression to obtain
\begin{eqnarray}\label{BanachContraction}
\|\theta_{n+1}\|_{L^\infty (L^p)}&\leq &\|e^{-\varepsilon t \mathcal{J}}\theta_0\|_{L^\infty (L^p)}+\left\|\int_{0}^t e^{-\varepsilon (t-s) \mathcal{J}}\nabla \cdot(v_\varepsilon\; \theta_n)(\cdot,s)ds \right\|_{L^\infty (L^p)}\nonumber\\
&& +\left\|\int_{0}^t e^{-\varepsilon (t-s) \mathcal{J}}\mathcal{J}^{1/2}\theta_n(\cdot,s)ds \right\|_{L^\infty (L^p)}
\end{eqnarray}
For the first term above we note that, since $e^{-\varepsilon t  \mathcal{J}}$ is a contraction operator, the estimate $\|e^{-\varepsilon t  \mathcal{J}}\theta_0\|_{L^p}\leq \| \theta_0\|_{L^p}$ is valid for all function $\theta_0\in L^{p}(\mathbb{G})$ with $1\leq p\leq +\infty$, for all $t>0$ and all $\varepsilon>0$. Thus, we have
\begin{equation}\label{Maj1}
\|e^{-\varepsilon t  \mathcal{J}}\theta_0\|_{L^\infty (L^p)}\leq \|\theta_0\|_{L^p}.
\end{equation}
For the second term of (\ref{BanachContraction}) we have the following fact: if $\theta_n\in L^{\infty}([0,T']; L^{p}(\mathbb{G}))$ and if $v\in  L^{\infty}([0,T']; L^{\infty}(\mathbb{G}))$, then 
\begin{eqnarray}
\left\|\int_{0}^t e^{-\varepsilon (t-s) \mathcal{J}}\nabla \cdot(v_\varepsilon\; \theta_n)(\cdot,s)ds \right\|_{L^\infty (L^p)}& = &\underset{0<t<T'}{\sup} \left\|\int_{0}^t \nabla \cdot(v_\varepsilon \theta_n)\ast h_{\varepsilon (t-s)}(\cdot,s)ds\right\|_{L^p}\nonumber\\
&\leq & \underset{0<t<T'}{\sup}\int_{0}^t  \left\|v_\varepsilon \theta_n(\cdot,s)\right\|_{L^p} \left\|\nabla h_{\varepsilon (t-s)}\right\|_{L^1} ds\nonumber\\
&\leq & \underset{0<t<T'}{\sup}\int_{0}^t  \left\|v_\varepsilon(\cdot,s)\right\|_{L^\infty}  \left\| \theta_n(\cdot,s)\right\|_{L^p} C(\varepsilon(t-s))^{-1/2} ds\nonumber\\
&\leq &  \|v\|_{L^\infty (L^\infty)}\left\|f\right\|_{L^\infty (L^p)}  \underset{0<t<T'}{\sup}\int_{0}^t C(\varepsilon(t-s))^{-1/2} ds\nonumber\\
&\leq &C \sqrt{\frac{T'}{\varepsilon}}  \|v\|_{L^\infty (L^\infty)} \left\| \theta_n\right\|_{L^\infty (L^p)}.\label{Maj3}
\end{eqnarray}
For the last term of (\ref{BanachContraction}) if $\theta_n\in L^{\infty}([0,T']; L^{p}(\mathbb{G}))$, then
\begin{eqnarray}
\left\|\int_{0}^t e^{-\varepsilon (t-s) \mathcal{J}}\mathcal{J}^{1/2}\theta_n(\cdot,s)ds \right\|_{L^\infty (L^p)}&=&\underset{0<t<T'}{\sup} \left\|\int_{0}^t \mathcal{J}^{1/2}  \theta_n\ast h_{\varepsilon (t-s)}(\cdot,s)ds\right\|_{L^p}\nonumber\\
&\leq & \underset{0<t<T'}{\sup} \displaystyle{\int_{0}^t} \|\theta_n(\cdot, s)\|_{L^p} \|\mathcal{J}^{1/2}h_{\varepsilon(t-s)}\|_{L^1}ds\nonumber\\
&\leq & C \left(\frac{T'^{1/2}}{\varepsilon^{1/2}}\right)\; \|\theta_n\|_{L^\infty (L^p)}\label{Maj2}
\end{eqnarray}

Now, applying the inequalities (\ref{Maj1}), (\ref{Maj3}) and (\ref{Maj2}) to the left-hand side of (\ref{BanachContraction}) we have
$$ \|\theta_{n+1}\|_{L^\infty (L^p)}\leq \|\theta_0\|_{L^p}+C\left(\frac{T'^{1/2}}{\varepsilon^{1/2}}\|v\|_{L^\infty (L^\infty)}+\frac{T'^{1/2}}{\varepsilon^{1/2}} \right)\|\theta_n\|_{L^\infty (L^p)}$$
Thus, if $\|\theta_0\|_{L^p}\leq K$ and if we define the time $T'$ to be such that $C\left(\frac{T'^{1/2}}{\varepsilon^{1/2}}+\frac{T'^{1/2}}{\varepsilon^{1/2}}\|v\|_{L^\infty (L^\infty)} \right)\leq 1/2$, we have by iteration that $\|\theta_{n+1}\|_{L^\infty (L^p)}\leq 2 K$: the sequence $(\theta_n)_{n\in \mathbb{N}}$ constructed from initial data $\theta_0$ belongs to the closed ball $\overline{B}(0, 2K)$. In order to finish this proof, let us show that $\theta_n \longrightarrow \theta$ in $L^{\infty}([0,T']; L^{p}(\mathbb{G}))$. For this we write
$$\|\theta_{n+1}-\theta_n\|_{L^\infty (L^p)}\leq \left\|\int_{0}^t e^{-\varepsilon (t-s) \mathcal{J}}\nabla \cdot(v_\varepsilon\; (\theta_{n}-\theta_{n-1}))(\cdot,s)ds\right\|_{L^\infty (L^p)}+\left\|\int_{0}^t e^{-\varepsilon (t-s) \mathcal{J}}\mathcal{J}^{1/2}(\theta_{n}-\theta_{n-1})\right\|_{L^\infty (L^p)}$$
and using the previous results we have
$$\|\theta_{n+1}-\theta_n\|_{L^\infty (L^p)}\leq C\left(\frac{T'^{1/2}}{\varepsilon^{1/2}}\|v\|_{L^\infty (L^\infty)}+\frac{T'^{1/2}}{\varepsilon^{1/2}}\right)\|\theta_{n}-\theta_{n-1}\|_{L^\infty (L^p)}$$
so, by iteration we obtain
$$\|\theta_{n+1}-\theta_n\|_{L^\infty (L^p)}\leq \left[C\left(\frac{T'^{1/2}}{\varepsilon^{1/2}}\|v\|_{L^\infty (L^\infty)}+\frac{T'^{1/2}}{\varepsilon^{1/2}}\right)\right]^{n}\|\theta_{1}-\theta_0\|_{L^\infty (L^p)}$$
hence, with the definition of $T'$ it comes
$\|\theta_{n+1}-\theta_n\|_{L^\infty (L^p)}\leq \left(\frac{1}{2}\right)^{n}\|\theta_{1}-\theta_0\|_{L^\infty (L^p)}$.
Finally, if $n\longrightarrow +\infty$, the sequence $(\theta_n)_{n\in \mathbb{N}}$ convergences towards $\theta$ in $L^\infty([0,T'];L^p(\mathbb{G}))$. Since it is a Banach space we deduce uniqueness for the solution $\theta$ of problem (\ref{FormIntegr}). The proof of Theorem \ref{TheoPointFixe} is finished.\hfill$\blacksquare$
\begin{Remarque}[From Local to Global]
Once we obtain a local result, global existence easily follows by a simple iteration since problems studied here (equations (\ref{Equation0}) or (\ref{SistApprox})) are linear as the velocity $v$ does not depend on $\theta$.
\end{Remarque}
We study now the regularity of the solutions constructed by this method.
\begin{Theoreme} Solutions of the approximated problem (\ref{SistApprox}) are smooth.
\end{Theoreme}
\textit{\textbf{Proof}.}
By iteration we will prove that $\theta \in \displaystyle{\bigcap_{0<T_0<T_1<t<T_2<T^\ast}}L^\infty([0,t]; W^{\frac{k}{2},p}(\mathbb{G}))$ for all $k\geq 0$. Note that this is true for $k=0$. So let us assume that it is also true for $k>0$ and we will show that it is still true for $k+1$.\\

Set $t$ such that $0<T_0<T_1<t<T_2<T^\ast$ and let us consider the next problem
$$\theta(x,t)=e^{-\varepsilon (t-T_0)\mathcal{J}}\theta(x, T_0)-\int_{T_0}^t e^{-\varepsilon (t-s)\mathcal{J}}\nabla \cdot(v_\varepsilon\; \theta)(x,s)ds-\int_{T_0}^t e^{-\varepsilon (t-s)\mathcal{J}}\mathcal{J}^{1/2} \theta(x,s)ds
$$
We have then the following estimate
\begin{eqnarray*}
\|\theta\|_{L^\infty (W^{\frac{k+1}{2},p})}&\leq& \|e^{-\varepsilon (t-T_0)\mathcal{J}}\theta(\cdot, T_0)\|_{L^\infty (W^{\frac{k+1}{2},p})}\\[5mm]
& &+\left\|\int_{T_0}^t e^{-\varepsilon (t-s)\mathcal{J}}\nabla \cdot(v_\varepsilon\; \theta)(\cdot,s)ds\right\|_{L^\infty (W^{\frac{k+1}{2},p})}
+\left\|\int_{T_0}^t e^{-\varepsilon (t-s)\mathcal{J}}\mathcal{J}^{1/2}\theta(\cdot,s)ds\right\|_{L^\infty (W^{\frac{k+1}{2},p})}
\end{eqnarray*}
Now, we will treat separately each of the previous terms. 
\begin{enumerate}
\item[(i)] For the first one we have
\begin{eqnarray*}
\|e^{-\varepsilon (t-T_0)\mathcal{J}}\theta(\cdot, T_0)\|_{W^{\frac{k+1}{2},p}}&=&\|\theta(\cdot, T_0)\ast h_{\varepsilon (t-T_0)}\|_{L^p}+\|\theta(\cdot, T_0)\ast\mathcal{J}^{\frac{k+1}{4}}h_{\varepsilon (t-T_0)}\|_{L^p}\\
&\leq &\|\theta(\cdot, T_0)\|_{L^p}+\|\theta(\cdot, T_0)\|_{L^p}\|\mathcal{J}^{\frac{k+1}{4}}h_{\varepsilon (t-T_0)}\|_{L^1} 
\end{eqnarray*}
where $h_t$ is the heat kernel, so we can write
\begin{equation*}
\|e^{-\varepsilon (t-T_0)\mathcal{J}}\theta(\cdot, T_0)\|_{L^\infty (W^{\frac{k+1}{2},p})}\leq C\|\theta(\cdot, T_0)\|_{L^p}\underset{T_1<t<T_2}{\sup}\left\{ \left[\varepsilon (t-T_0)\right]^{- \frac{k+1}{4}}; 1\right\}
\end{equation*}
\item[(ii)] For the second term, one has
\begin{eqnarray*}
I&=&\left\|\int_{T_0}^t e^{-\varepsilon (t-s)\mathcal{J}}\nabla \cdot(v_\varepsilon\; \theta)(\cdot,s)ds\right\|_{W^{\frac{k+1}{2},p}}\leq \int_{T_0}^t\|\nabla  \cdot (v_{\varepsilon}\;\theta)\ast h_{\varepsilon(t-s)}\|_{L^p}+\|\nabla  \cdot (v_{\varepsilon}\;\theta)\ast h_{\varepsilon(t-s)}\|_{\dot{W}^{\frac{k+1}{2},p}}ds\\
&\leq &\int_{T_0}^t \|v_{\varepsilon}\;\theta\|_{L^p}\|\nabla h_{\varepsilon(t-s)}\|_{L^1}+ \|\mathcal{J}^{\frac{k}{4}}(v_{\varepsilon}\;\theta)\|_{L^p}\|\mathcal{J}^{1/4}\big(\nabla h_{\varepsilon(t-s)}\big)\|_{L^1}ds\\
&\leq &C\int_{T_0}^t\|v_\varepsilon\;\theta(\cdot,s)\|_{L^p} \left[\varepsilon (t-s)\right]^{- \frac{1}{2}}+\|v_\varepsilon\;\theta(\cdot,s)\|_{\dot{W}^{\frac{k}{2},p}} \left[\varepsilon (t-s)\right]^{- \frac{3}{4}} ds.\\
&\leq &C\int_{T_0}^t\|v_\varepsilon\;\theta(\cdot,s)\|_{W^{\frac{k}{2},p}}\max\left(\left[\varepsilon (t-s)\right]^{- \frac{1}{2}}; \left[\varepsilon (t-s)\right]^{- \frac{3}{4}}\right)ds
\end{eqnarray*}
Note now that we have here the estimations below for $k/2\leq \ell\in \mathbb{N}$
\begin{eqnarray*}
\|v_\varepsilon\theta(\cdot,s)\|_{W^{\frac{k}{2},p}}&\leq& \|v_\varepsilon(\cdot,s)\|_{\mathcal{C}^\ell} \|\theta(\cdot,s)\|_{W^{\frac{k}{2},p}}\leq C \varepsilon^{-\ell}\|v(\cdot,s)\|_{L^\infty}\|\theta(\cdot,s)\|_{W^{\frac{k}{2},p}}
\end{eqnarray*}
hence, we can write
\end{enumerate}
\begin{eqnarray*}
\left\|\int_{T_0}^t e^{-\varepsilon (t-s)\mathcal{J}}\nabla \cdot(v_\varepsilon\; \theta)(\cdot,s)ds\right\|_{L^\infty (W^{\frac{k+1}{2},p})}\leq C \|v\|_{L^\infty (L^\infty)} \|\theta\|_{L^\infty (W^{\frac{k}{2},p})}\underset{T_1<t<T_2}{\sup}\int_{T_0}^t  \varepsilon^{-\ell}\max\left(\left[\varepsilon (t-s)\right]^{- \frac{1}{2}}; \left[\varepsilon (t-s)\right]^{- \frac{3}{4}} \right)ds
\end{eqnarray*}
\begin{enumerate}
\item[(iii)] Finally, for the last term we have
\begin{eqnarray*}
\left\|\int_{T_0}^t e^{-\varepsilon (t-s)\mathcal{J}}\mathcal{J}^{1/2}\theta(\cdot,s)ds\right\|_{W^{\frac{k+1}{2},p}} &\leq & \int_{T_0}^t \|\theta(\cdot,s)\|_{L^p}\|\mathcal{J}^{1/2}h_{\varepsilon(t-s)}\|_{L^1}+ \left\|\mathcal{J}^{\frac{k}{4}}\theta(\cdot,s)\right\|_{L^{p}}\|\mathcal{J}^{3/4}h_{\varepsilon(t-s)}\|_{L^{1}}ds\\
&\leq & C  \int_{T_0}^t \|\theta(\cdot,s)\|_{W^{\frac{k}{2},p}}\max\left (\left[\varepsilon (t-s)\right]^{- \frac{1}{2}}; \left[\varepsilon (t-s)\right]^{- \frac{3}{4}}\right)ds
\end{eqnarray*}
So finally we have
\begin{eqnarray*}
\left\|\int_{T_0}^t e^{-\varepsilon (t-s)\mathcal{J}} \mathcal{J}^{1/2}  \theta(\cdot,s)ds\right\|_{L^\infty (W^{\frac{k+1}{2},p})}&\leq & C\|\theta \|_{L^{\infty}(W^{\frac{k}{2},p})}\underset{T_1<t<T_2}{\sup}\int_{T_0}^t\max\left (\left[\varepsilon (t-s)\right]^{- \frac{1}{2}}; \left[\varepsilon (t-s)\right]^{- \frac{3}{4}}\right)ds.
\end{eqnarray*}
\end{enumerate}
Now, with formulas (i)-(iii) at our disposal, we have that the norm $\|\theta\|_{L^\infty (W^{\frac{k+1}{2},p})}$ is controlled for all $\varepsilon>0$: we have proven spatial regularity.\\

Time regularity follows since we have
$$\frac{ \partial^k}{\partial t^k}\theta(x,t)+\nabla \cdot \left(\frac{ \partial^k}{\partial t^k} (v_\varepsilon\,\theta)\right)(x,t)+\mathcal{J}^{1/2}\left(\frac{ \partial^k}{\partial t^k} \theta\right)(x,t)=\varepsilon \Delta \left(\frac{ \partial^k}{\partial t^k} \theta\right)(x,t).$$
\begin{flushright}$\blacksquare$\end{flushright}
\begin{Remarque}\label{RemarkEpsilonDep}
The solutions $\theta(\cdot,\cdot)$ constructed above depend on $\varepsilon$. 
\end{Remarque}
\subsection{Maximum principle for regular solutions}\label{SeccMax} 
The maximum principle we are studying here will be a consequence of few inequalities. We will start with the solutions obtained in the previous section:
\begin{Proposition}[Regularized version of Theorem \ref{Theo1}-1)]\label{PropoViscosityMaxPrinc}
Let $\theta_0\in L^{p}(\mathbb{G})$ with $1\leq p\leq +\infty$ be an initial data, then the associated solution of the problem  (\ref{SistApprox}) satisfies the maximum principle for all $t\in [0, T]$: $\|\theta(\cdot, t)\|_{L^p}\leq \|\theta_0\|_{L^p}$.
\end{Proposition}
\textit{\textbf{Proof.}} We write for $1\leq p<+\infty$:
\begin{eqnarray*}
\frac{d}{dt}\|\theta(\cdot, t)\|^p_{L^p}&=&p\int_{\mathbb{G}}|\theta|^{p-2}\theta\bigg(-\varepsilon \mathcal{J} \theta-\nabla \cdot (v_\varepsilon \,\theta)-\mathcal{J}^{1/2}\theta\bigg)dx=-p\varepsilon\int_{\mathbb{G}}|\theta|^{p-2}\theta\mathcal{J} \theta dx-p\int_{\mathbb{G}}|\theta|^{p-2}\theta \mathcal{J}^{1/2}\theta dx\\
\end{eqnarray*}
where we used the fact that $div(v)=0$. Thus, we have
$$\frac{d}{dt}\|\theta(\cdot, t)\|^p_{L^p}+p\varepsilon\int_{\mathbb{G}}|\theta|^{p-2}\theta \mathcal{J} \theta dx+p\int_{\mathbb{G}}|\theta|^{p-2}\theta\mathcal{J}^{1/2}\theta dx=0,$$
and integrating in time we obtain
\begin{equation}\label{Form1}
\|\theta(\cdot, t)\|^p_{L^p}+p\varepsilon\int_{0}^t\int_{\mathbb{G}}|\theta|^{p-2}\theta \mathcal{J} \theta dxds+p\int_0^t\int_{\mathbb{G}}|\theta|^{p-2}\theta\mathcal{J}^{1/2} \theta dxds=\|\theta_0\|^p_{L^p}.
\end{equation}
To finish, we have that the quantities $p\varepsilon\displaystyle{\int_{\mathbb{G}}}|\theta|^{p-2}\theta \mathcal{J} \theta dx\quad \mbox{ and }\quad \displaystyle{\int_0^t\int_{\mathbb{G}}}|\theta|^{p-2}\theta \mathcal{J}^{1/2} \theta dxds$ are both positive. For the first expression, this is a consecuence of the fact that $e^{-\varepsilon s \mathcal{J}}$  is a contraction semi-group. For the second expression we can use the Positivity Lemma of \cite{Cordoba} which is valid for $1\leq p<+\infty$; however, we will need the following lemma (see a proof in \cite{PGDCH}):
\begin{Lemme}\label{TheoBesov}
If $2\leq p<+\infty$, then there is a positive constant $C>0$ such that
\begin{equation*}
C \|\theta\|_{\dot{B}^{1/p, p}_p}^p\leq \int_{\mathbb{G}}|\theta|^{p-2}\theta \mathcal{J}^{1/2} \theta dx
\end{equation*}
\end{Lemme}
Thus, getting back to (\ref{Form1}), we have that all these quantities are bounded and positive and we write for all $1\leq p<+\infty$:
$$\|\theta(\cdot, t)\|_{L^p}\leq \|\theta_0\|_{L^p}.$$
Since $\|\theta(\cdot, t)\|_{L^p}\underset{p\to +\infty}{\longrightarrow}\|\theta(\cdot, t)\|_{L^\infty}$, the maximum principle is proven for regular solutions.  \hfill$\blacksquare$\\

Of course, this remains true for smooth solutions of equation (\ref{Equation0}). 

\subsection{The limit $\varepsilon\longrightarrow 0$ for regular solutions} 
We have proven so far regular versions (i.e. for $\varepsilon>0$) of Theorem \ref{Theo0} and Theorem \ref{Theo1}. We will now pass to the limit  $\varepsilon\longrightarrow 0$.\\

\textit{\textbf{Proof of Theorem \ref{Theo0}.}} We have obtained with the previous results a family of regular functions $(\theta^{(\varepsilon)})_{\varepsilon >0}\in L^{\infty}([0,T]; L^{p}(\mathbb{G}))$ which are solutions of (\ref{SistApprox}) and satisfy the uniform bound $\|\theta^{(\varepsilon)}(\cdot, t)\|_{L^p}\leq \|\theta_0\|_{L^p}$. \\

Since $L^{\infty}([0,T]; L^{p}(\mathbb{G}))= \left(L^{1}([0,T]; L^{q}(\mathbb{G}))\right)'$, with $\frac{1}{p}+\frac{1}{q}=1$, we can extract from those solutions $\theta^{(\varepsilon)}$ a subsequence $(\theta_k)_{k\in \mathbb{N}}$ which is $\ast$-weakly convergent to some function $\theta$ in the space $L^{\infty}([0,T]; L^{p}(\mathbb{G}))$, which implies convergence in $\mathcal{D}'(\mathbb{R}^+\times\mathbb{G})$. However, this weak convergence is not sufficient to assure the convergence of $(v_\varepsilon\; \theta_k)$ to $v\;\theta$. For this we use the  remarks that follow.\\

First, using the Lemma \ref{TheoApproxbmo} we can consider a sequence $(v_k)_{k\in\mathbb{N}}$ such that $v_k \longrightarrow v$ weakly in $bmo(\mathbb{G})$. Secondly, combining  Proposition \ref{PropoViscosityMaxPrinc} and Lemma \ref{TheoBesov} we obtain that solutions $\theta_k$ belongs to the space $L^{\infty}([0,T]; L^{p}(\mathbb{G}))\cap L^{1}([0,T]; \dot{B}^{1/p,p}_p(\mathbb{G}))$ for all $k\in \mathbb{N}$. \\

To finish, fix a function $\varphi\in \mathcal{C}^{\infty}_{0}([0,T]\times \mathbb{G})$. Then we have the fact that $\varphi \theta_k\in  L^{1}([0,T]; \dot{B}^{1/p,p}_p(\mathbb{G}))$ and $\partial_t \varphi \theta_k\in  L^{1}([0,T]; \dot{B}^{-\ell,p}_p(\mathbb{G}))$. This implies the local inclusion, in space as well as in time, $\varphi \theta_k\in \dot{W}^{1/p,p}_{t,x}\subset \dot{W}^{1/p,2}_{t,x}$ so we can apply classical results such as the Rellich's theorem to obtain convergence of $v_k\; \theta_k$ to $v\;\theta$. \\

Thus, we obtain existence and uniqueness of weak solutions for the problem (\ref{Equation0}) with an initial data in $\theta_0\in L^p(\mathbb{G})$, $2\leq p<+\infty$ that satisfy the maximum principle. Moreover, we have that these solutions $\theta(x,t)$ belong to the space $L^{\infty}([0,T]; L^{p}(\mathbb{G}))\cap L^{p}([0,T]; \dot{B}^{1/p,p}_p(\mathbb{G}))$.\hfill$\blacksquare$
\begin{Remarque}
These lines explain how to obtain weak solutions from viscosity ones and it will be used freely in the sequel.  
\end{Remarque}
\section{Positivity principle}\label{Sect_PrincipeMax}
We prove in this section Theorem \ref{Theo1}-2). Recall that we have $0\leq \psi_0\leq M$ with $\psi_0\in L^p(\mathbb{G})$ and $N< p\leq+\infty$, where $N$ is the homogeneous dimension of $\mathbb{G}$. We will show that the associated solution $\psi(x,t)$ satisfies the bounds $0\leq \psi(x,t)\leq M$.\\

To begin with, we fix two constants, $\rho, R$ such that $R>2\rho>0$. Then we set $A_{0,R}(x)$ a function equals to $M/2$ over $\|x\|\leq 2R$ and equals to $\psi_0(x)$ over $\|x\|>2R$ and we write $B_{0,R}(x)=\psi_0(x)-A_{0,R}(x)$, so by construction we have $$\psi_0(x)=A_{0,R}(x)+B_{0,R}(x)$$ with $\|A_{0,R}\|_{L^\infty} \leq M$ and $\|B_{0,R}\|_{L^\infty} \leq M/2$. Remark that $A_{0,R}, B_{0,R}\in L^p(\mathbb{G})$.\\

Now fix $v\in L^{\infty}([0,T];bmo(\mathbb{G}))$ such that $div(v)=0$ and consider the equations
\begin{equation}\label{ForDouble}
\begin{cases}
\partial_t A_R+ \nabla\cdot (v\,A_R)+\mathcal{J}^{1/2} A_R=0,\\[4mm]
A_R(x,0)=A_{0,R}(x)
\end{cases}
\mbox{and}\qquad
\begin{cases}
\partial_t B_R+ \nabla\cdot (v\,B_R)+\mathcal{J}^{1/2} B_R=0\\[4mm]
B_R(x,0)= B_{0,R}(x).
\end{cases}
\end{equation}
Using the maximum principle and by construction we have the following estimates for $t\in [0,T]$:
\begin{eqnarray}
\|A_R(\cdot, t)\|_{L^p}&\leq& \|A_{0,R}\|_{L^p}\leq \|\psi_0\|_{L^p}+CM R^{N/p} \quad (1<p<+\infty)\label{Formula1}\\[5mm] 
\|A_R(\cdot, t)\|_{L^\infty}&\leq & \|A_{0,R}\|_{L^\infty}\leq M.\nonumber\\[5mm]
\|B_R(\cdot, t)\|_{L^\infty}&\leq & \|B_{0,R}\|_{L^\infty}\leq M/2. \nonumber
\end{eqnarray}
where $A_R(x,t)$ and $B_R(x,t)$ are the weak solutions of the systems (\ref{ForDouble}). Then, the function $\psi(x,t)=A_R(x,t)+B_R(x,t)$ is the unique solution for the problem
\begin{equation}\label{Equation1}
\left\lbrace
\begin{array}{l}
\partial_t \psi+ \nabla\cdot (v\,\psi)+\mathcal{J}^{1/2} \psi=0\\[5mm]
\psi(x,0)=A_{0,R}(x)+B_{0,R}(x).
\end{array}
\right.
\end{equation}
Indeed, using hypothesis for $A_R(x,t)$ and $B_R(x,t)$ and the linearity of equation (\ref{Equation1}) we have that the function $\psi_R(x,t)=A_R(x,t)+B_R(x,t)$ is a solution for this equation. Uniqueness is assured by the maximum principle, thus we can write $\psi(x,t)=\psi_R(x,t)$.\\

To continue, we will need an auxiliary function $\phi \in \mathcal{C}^{\infty}_{0}(\mathbb{G})$ such that $\phi(x)=0$ for $\|x\|\geq 1$ and $\phi(x)=1$ if $\|x\|\leq 1/2$ and we set $\varphi(x)=\phi(\delta_{R^{-1}}[x])$. Now, we will estimate the $L^p$-norm of $\varphi(x)(A_R(x,t)-\frac{M}{2})$ with $p>N$. We write:
\begin{eqnarray}
\partial_t\left\|\varphi(\cdot)\left(A_R(\cdot,t)-\frac{M}{2}\right)\right\|_{L^p}^p &=& p\int_{\mathbb{G}}\left|\varphi(x)(A_R(x,t)-\frac{M}{2})\right|^{p-2}\big(\varphi(x)(A_R(x,t)-\frac{M}{2}) \big)\nonumber\\[4mm]
& & \times\; \partial_t\big(\varphi(x)(A_R(x,t)-\frac{M}{2}) \big) dx \label{equat1}
\end{eqnarray}
We observe that we have the following identity for the last term above
\begin{eqnarray*}
\partial_t\left(\varphi(x)\left(A_R(x,t)-\frac{M}{2}\right)\right)&=&- \nabla \cdot \left(\varphi(x) \, v(A_R(x,t)-\frac{M}{2})\right)-\mathcal{J}^{1/2} \left(\varphi(x)(A_R(x,t)-\frac{M}{2})\right)\\[3mm]
&+&\left(A_R(x,t)-\frac{M}{2}\right)v\cdot \nabla \varphi(x)+[\mathcal{J}^{1/2}, \varphi]A_R(x,t)-\frac{M}{2} \mathcal{J}^{1/2} \varphi(x)
\end{eqnarray*}
where we noted $[\mathcal{J}^{1/2}, \varphi]$ the commutator between $\mathcal{J}^{1/2}$ and $\varphi$. Thus, using this identity in (\ref{equat1}) and the fact that $div(v)=0$ we have
\begin{eqnarray}
\partial_t\left\|\varphi(\cdot)\left(A_R(\cdot,t)-\frac{M}{2}\right)\right\|_{L^p}^p &=&-p\int_{\mathbb{G}}\big|\varphi(x)(A_R(x,t)-\frac{M}{2})\big|^{p-2}\big(\varphi(x)(A_R(x,t)-\frac{M}{2}) \big)\nonumber\\[4mm]
& & \times\; \mathcal{J}^{1/2}(\varphi(x)\big(A_R(x,t)-\frac{M}{2})\big)dx\label{equa2}\\[4mm]
&+&p\int_{\mathbb{G}}\big|\varphi(x)(A_R(x,t)-\frac{M}{2})\big|^{p-2}\big(\varphi(x)(A_R(x,t)-\frac{M}{2}) \big)\nonumber\\[4mm]
& & \times\; \left([\mathcal{J}^{1/2}, \varphi]A_R(x,t)-\frac{M}{2}\mathcal{J}^{1/2} \varphi(x)\right)dx\nonumber
\end{eqnarray}
Remark that the integral (\ref{equa2}) is positive so one has
\begin{eqnarray*}
\partial_t\left\|\varphi(\cdot)\left(A_R(\cdot,t)-\frac{M}{2}\right)\right\|_{L^p}^p &\leq &p\int_{\mathbb{G}}\big|\varphi(x)(A_R(x,t)-\frac{M}{2})\big|^{p-2}\big(\varphi(x)(A_R(x,t)-\frac{M}{2}) \big)\\[4mm]
& & \times\; \left([\mathcal{J}^{1/2}, \varphi]A_R(x,t)-\frac{M}{2} \mathcal{J}^{1/2}\varphi(x)\right)dx
\end{eqnarray*}
Using Hölder's inequality and integrating in time the previous expression we have
\begin{eqnarray}\label{EquaPrevious}
\left\|\varphi(\cdot)\left(A_R(\cdot,t)-\frac{M}{2}\right)\right\|^p_{L^p}& \leq &\left\|\varphi(\cdot)\big(A_R(\cdot,0)-\frac{M}{2}\big)\right\|^p_{L^p}+\int_{0}^t\bigg(\left\|[\mathcal{J}^{1/2}, \varphi]A_R(\cdot,s)\right\|_{L^p}+ \|\frac{M}{2} \mathcal{J}^{1/2}\varphi\|_{L^p}\bigg)ds
\end{eqnarray}
The first term of the right side is null since over the support of $\varphi$ we have identity $A_R(x,0)=\frac{M}{2}$. For the term $\left\|[\mathcal{J}^{1/2}, \varphi]A_R(\cdot,s)\right\|_{L^p}$ we will need the following lemma (see a proof in \cite{Stein2} or \cite{Grafakos}):
\begin{Lemme}\label{LemmaInterpolCom}
For $1< p\leq +\infty$ we have the following inequality:
\begin{equation*}
\left\|[\mathcal{J}^{1/2}, \varphi]A_R(\cdot,s)\right\|_{L^p}\leq CR^{-1}\|A_{0,R}\|_{L^p}.
\end{equation*}
\end{Lemme}
Now, getting back to the last term of (\ref{EquaPrevious}) we have by the definition of $\varphi$ and the properties of the operator $\mathcal{J}^{1/2}$ the estimate:
$$\|\frac{M}{2} \mathcal{J}^{1/2}\varphi\|_{L^p}\leq C MR^{N/p}R^{-1}.$$
We thus have
$$\left\|\varphi(\cdot)\left(A_R(\cdot,t)-\frac{M}{2}\right)\right\|^p_{L^p}\leq CR^{-1}\int_{0}^t\bigg(\|A_{0,R}\|_{L^p}+MR^{N/p}\bigg)ds.$$
Observe that we have at our disposal estimate (\ref{Formula1}), so we can write
$$\left\|\varphi(\cdot)\left(A_R(\cdot,t)-\frac{M}{2}\right)\right\|^p_{L^p}\leq Ct R^{-1}\left(\|\psi_0\|_{L^p}+MR^{N/p}\right)$$
Using again the definition of $\varphi$ one has $\displaystyle{\int_{B(0,\rho)}}|A_R(\cdot,t)-\frac{M}{2}|^{p}dx\leq CtR^{-1}\left(\|\psi_0\|_{L^p}+MR^{N/p}\right)$. Thus, if $R\longrightarrow +\infty$ and since $p>N$, we have $A(x,t)=\frac{M}{2}$ over $B(0,\rho)$.\\ 

Hence, by construction we have  $\psi(x,t)=A_R(x,t)+B_R(x,t)$ where $\psi$ is a solution of (\ref{Equation1}) with initial data $\psi_0=A_{0,R}+B_{0,R}$, but, since over $B(0,\rho)$ we have $A(x,t)=\frac{M}{2}$ and $\|B(\cdot,t)\|_{L^\infty}\leq \frac{M}{2}$, one finally has the desired estimate $0\leq \psi(x,t)\leq M$.\hfill$\blacksquare$
\section{Existence of solutions with a $L^\infty$ initial data}\label{SecLinfty}

The proof given before for the positivity principle allows us to obtain the existence of solutions for the fractional diffusion transport equation (\ref{Equation0}) when the initial data $\theta_0$ belongs to the space $L^\infty(\mathbb{G})$. The utility of this fact will appear clearly in the next section as it will be used in Theorem \ref{Theo3}.\\

Let us fix $\theta_0^R=\theta_0 \mathds{1}_{B(0,R)}$ with $R>0$ so we have $\theta_0^R\in L^p(\mathbb{G})$ for all $1\leq p\leq +\infty$. Following section \ref{SecExiUnic}, there is a unique solution $\theta^R$ for the problem
\begin{equation*}
\left\lbrace
\begin{array}{l}
\partial_t \theta^R+ \nabla\cdot (v\theta^R)+\mathcal{J}^{1/2}\theta^R=0\\[5mm]
\theta^R(x,0)=\theta_0^R(x)\\[5mm]
div(v)=0 \quad \mbox{ and } v\in L^{\infty}([0,T];  bmo(\mathbb{G})).
\end{array}
\right.
\end{equation*}
such that $\theta^R\in L^\infty([0,T]; L^p(\mathbb{G}))$. By the maximum principle we have $\|\theta^R(\cdot, t)\|_{L^p}\leq\|\theta^R_0\|_{L^p}\leq v_n\|\theta_0\|_{L^\infty}R^{N/p}$. Taking the limit $p\longrightarrow +\infty$ and making $R\longrightarrow +\infty$ we finally get
$$\|\theta(\cdot, t)\|_{L^\infty}\leq C \|\theta_0\|_{L^\infty}.$$
This shows that for an initial data $\theta_0\in  L^\infty(\mathbb{G})$ there exists an associated solution $\theta\in L^\infty([0,T];L^\infty(\mathbb{G}))$.
\section{H\"older Regularity}\label{SeccHolderRegularity}

In this section we are going to prove Theorem \ref{Theo3}. Our aim is to prove that the solutions of equation (\ref{Equation0}) are $\gamma$-H\"older regular (with $0<\gamma<1/2 $) in the sense that the following norm is bounded
$$\|\theta\|_{\mathcal{C}^\gamma}=\|\theta\|_{L^\infty}+\underset{x,y\in \mathbb{G}}{\sup}\frac{|\theta(x\cdot y)- \theta(x)|}{\|y\|^\gamma}$$
However, we will not work with this quantity, we will use instead a duality characterization based on Hardy spaces $h^\sigma$ with  $0<\sigma<1$. Indeed in the framework of stratified Lie groups we have that $(h^\sigma)'=\mathcal{C}^\gamma$ with $\gamma=N(\frac{1}{\sigma}-1)$ (see \cite{Folland2} for a proof), thus in order to verify that a function $\theta$ is H\"older regular it is enough to study the quantity
\begin{equation*}
\langle \theta, \psi \rangle_{\mathcal{C}^\gamma \times h^\sigma}
\end{equation*}
Hardy spaces $h^\sigma$ have several equivalent characterizations and in this paper we are interested mainly in the molecular approach that defines \textit{local} Hardy spaces.
\begin{Definition}[Local Hardy spaces $h^\sigma$] Let $0<\sigma<1$. The local Hardy space $h^{\sigma}(\mathbb{G})$ is the set of distributions $f$ that admits the following molecular decomposition:
\begin{equation}\label{MolDecomp}
f=\sum_{j\in \mathbb{N}}\lambda_j \psi_j
\end{equation}
where $(\lambda_j)_{j\in \mathbb{N}}$ is a sequence of complex numbers such that $\sum_{j\in \mathbb{N}}|\lambda_j|^\sigma<+\infty$ and $(\psi_j)_{j\in \mathbb{N}}$ is a family of $r$-molecules in the sense of definition \ref{DefMolecules} below. The $h^\sigma$-norm\footnote{it is not actually a \textit{norm} since $0<\sigma<1$. More details can be found in \cite{Gold} and \cite{Stein2}.} is then fixed by the formula 
$$\|f\|_{h^\sigma}=\inf\left\{\left(\sum_{j\in \mathbb{N}}|\lambda_j|^\sigma\right)^{1/\sigma}:\; f=\sum_{j\in \mathbb{N}}\lambda_j \psi_j \right\}$$ where the infimum runs over all possible decompositions (\ref{MolDecomp}).
\end{Definition}
Local Hardy spaces have many remarquable properties and we will only stress here, before passing to duality results concerning $h^\sigma$ spaces, the fact that Schwartz class $\mathcal{S}(\mathbb{G})$ is dense in $h^{\sigma}(\mathbb{G})$. 
\begin{Remarque}\label{Remark3}
Since $0<\sigma<1$, we have $\sum_{j\in \mathbb{N}}|\lambda_j|\leq\left(\sum_{j\in \mathbb{N}}|\lambda_j|^\sigma\right)^{1/\sigma}$ thus for testing Hölder continuity of a function $f$ it is enough to study the quantities $|\langle f,\psi_j\rangle|$ where $\psi_j$ is an $r$-molecule.
\end{Remarque}
Molecules in stratified Lie groups have been studied in \cite{DeMichele}. We give here an equivalent definition which is more suited to our purposes. 
\begin{Definition}[$r$-molecules]\label{DefMolecules} Set $\frac{N}{N+1}<\sigma<1$, define $\gamma=N(\frac{1}{\sigma}-1)$ and fix a real number $\omega$ such that $0<\gamma<\omega<1$. An integrable function $\psi$ is an $r$-molecule if we have
\begin{enumerate}
\item[$\bullet$]\underline{Small molecules $(0<r<1)$:}
\begin{eqnarray}
& & \int_{\mathbb{G}} |\psi(x)|\|x\cdot x_0^{-1}\|^{\omega}dx \leq  r^{\omega-\gamma}\mbox{, for } x_0\in \mathbb{G}\label{Hipo1} \;\qquad\qquad\qquad\mbox{(concentration condition)} \\[5mm]
& &\|\psi\|_{L^\infty}  \leq  \frac{1}{r^{N+\gamma}}\label{Hipo2} \qquad\qquad\qquad\qquad\qquad\qquad\qquad\qquad\qquad\, \mbox{(height condition)}  \\[5mm]
& &\int_{\mathbb{G}} \psi(x)dx=0\label{Hipo3}\qquad\qquad\qquad\qquad\qquad\qquad\qquad\qquad\qquad \mbox{(moment condition)} 
\end{eqnarray}
\item[$\bullet$] \underline{Big molecules $(1\leq r<+\infty)$:}\\

In this case we only require conditions (\ref{Hipo1}) and (\ref{Hipo2}) for the $r$-molecule $\psi$ while the moment condition (\ref{Hipo3}) is dropped.
\end{enumerate}
\end{Definition}
\begin{Remarque}\label{Remark2}
\begin{itemize}
\item[]
\item[1)] Note that the point $x_0\in \mathbb{G}$ can be considered as the ``center'' of the molecule.
\item[2)] Conditions (\ref{Hipo1}) and (\ref{Hipo2}) imply the estimate $\|\psi\|_{L^1}\leq C\, r^{-\gamma}$ thus every $r$-molecule belongs to $L^p(\mathbb{G})$ with $1<p<+\infty$. In particular we have
\begin{equation*}
\| \psi\|_{L^p} \leq C r^{-N(1-1/p)-\gamma}.
\end{equation*}
\end{itemize}
\end{Remarque}
The main interest of using molecules relies in the possibility of \textit{transfering} the regularity problem to the evolution of such molecules:
\begin{Proposition}[Transfer property]\label{Transfert} Let $\psi(x,s)$ be a solution of the backward problem
\begin{equation}
\left\lbrace
\begin{array}{rl}
\partial_s \psi(x,s)=& -\nabla\cdot [v(x,t-s)\psi(x,s)]-\mathcal{J}^{1/2}\psi(x,s)\label{Evolution1}\\[5mm]
\psi(x,0)=& \psi_0(x)\in L^1\cap L^\infty(\mathbb{G})\\[5mm]
div(v)=&0 \quad \mbox{and }\; v\in L^{\infty}([0,T];bmo(\mathbb{G}))
\end{array}
\right.
\end{equation}
If $\theta(x,t)$ is a solution of (\ref{Equation0}) with $\theta_0\in L^\infty(\mathbb{G})$ then we have the identity
\begin{equation*}
\int_{\mathbb{G}}\theta(x,t)\psi(x,0)dx=\int_{\mathbb{G}}\theta(x,0)\psi(x,t)dx.
\end{equation*}
\end{Proposition}
\textit{\textbf{Proof.}}
We first consider the expression
$$\partial_s\int_{\mathbb{G}}\theta(x,t-s)\psi(x,s)dx=\int_{\mathbb{G}}-\partial_s\theta(x,t-s)\psi(x,s)+\partial_s\psi(x,s)\theta(x,t-s)dx.$$
Using equations (\ref{Equation0}) and (\ref{Evolution1}) we obtain
\begin{eqnarray*}
\partial_s\int_{\mathbb{G}}\theta(x,t-s)\psi(x,s)dx&=&\int_{\mathbb{G}}- \nabla\cdot\left[(v(x,t-s)\theta(x,t-s)\right]\psi(x,s)+\mathcal{J}^{1/2}\theta(x,t-s)\psi(x,s) \\[5mm]
&-& \nabla\cdot\left[(v(x,t-s) \psi(x,s))\right]\theta(x,t-s)-\mathcal{J}^{1/2}\psi(x,s) \theta(x,t-s)dx.
\end{eqnarray*}
Now, using the fact that $v$ is divergence free and the symmetry of the operator $\mathcal{J}^{1/2}$ we have that the expression above is equal to zero, so the quantity
$$\int_{\mathbb{G}}\theta(x,t-s)\psi(x,s)dx$$
remains constant in time. We only have to set $s=0$ and $s=t$ to conclude. \hfill$\blacksquare$\\

This proposition says, that in order to control $\langle \theta(\cdot, t),\psi_0\rangle$, it is enough (and much simpler) to study the bracket $\langle \theta_0,\psi(\cdot, t)\rangle$. \\

\textbf{Proof of Theorem \ref{Theo3}.} 
Once we have the transfer property proven above, the proof of Theorem \ref{Theo3} is quite direct and it reduces to a $L^q$ estimate for molecules with $\frac{1}{p}+\frac{1}{q}=1$. Indeed, assume that for \textit{all} molecular initial data $\psi_0$ we have a $L^q$ control for $\psi(\cdot,t)$ a solution of (\ref{Evolution1}), then Theorem \ref{Theo3} follows easily: applying Proposition \ref{Transfert} with the fact that $\theta_0\in L^{p}(\mathbb{G})$ we have 
\begin{equation}\label{DualQuantity1}
|\langle \theta(\cdot, t), \psi_0\rangle|=\left|\int_{\mathbb{G}}\theta(x,t)\psi_0(x)dx\right|=\left|\int_{\mathbb{G}}\theta(x,0)\psi(x,t)dx\right|\leq \|\theta_0\|_{L^p}\|\psi(\cdot,t)\|_{L^q}<+\infty.
\end{equation}
From this, we obtain that $\theta(\cdot, t)$ belongs to the Hölder space $\mathcal{C}^\gamma(\mathbb{G})$.\\

Now we need to study the control of the $L^q$ norm of $\psi(\cdot, t)$ and we divide our proof in two steps following the molecule's size. For the initial big molecules, \textit{i.e.} if $r\geq 1$, the needed control is straightforward: apply the maximum principle and the remark \ref{Remark2}-2) above to obtain
$$\|\theta_0\|_{L^p}\|\psi(\cdot,t)\|_{L^q}\leq \|\theta_0\|_{L^p}\|\psi_0\|_{L^q}\leq C \frac{1}{r^{N(1-1/p)+\gamma}} \|\theta_0\|_{L^\infty},$$
but, since $r\geq 1$, we have that $|\langle \theta(\cdot, t), \psi_0\rangle|<+\infty$ for all \textit{big} molecules.\\

In order to finish the proof of this theorem, it  only remains to treat the $L^q$ control for \textit{small} molecules. This is the most complex part of the proof and it is treated in the following theorem:
\begin{Theoreme}\label{TheoL1control}
For all small $r$-molecules (\textit{i.e. }$0<r<1$), there exists a time $T_0>0$ such that we have the following control of the $L^q$-norm.
$$\|\psi(\cdot, t)\|_{L^q}\leq C T_0^{-N(1-1/q)-\gamma}\qquad (T_0<t<T),$$
with $0<\gamma<1$. 
\end{Theoreme}
Accepting for a while this result, we have then a good control over the quantity $\|\psi(\cdot, t)\|_{L^q}$  for all $0<r<1$ and getting back to (\ref{DualQuantity1}) we obtain that $|\langle \theta(\cdot, t), \psi_0\rangle|$ is always bounded for $T_0<t<T$ and for any molecule $\psi_0$: we have proven by a duality argument the Theorem \ref{Theo3}. \hfill $\blacksquare$\\

Let us now briefly explain the main steps of Theorem \ref{TheoL1control}. We need to construct a suitable control in time for the $L^q$-norm of the solutions $\psi(\cdot,t)$ of the backward problem (\ref{Evolution1}) where the inital data $\psi_0$ is a \textit{small} $r$-molecule. This will be achieved by iteration in two different steps: 
\begin{itemize}
\item The first step explains the molecules' deformation after a very small time $s_0>0$. We will thus obtain similar concentration and height conditions from wich, applying remark \ref{Remark2}-2), we will obtain a $L^q$ bound for small times. This will be done in section \ref{SecEvolMol1}.
\item In order to obtain a control of the $L^q$ norm for larger times we need to perform a second step which takes as a starting point the results of the first step and gives us the deformation for another small time $s_1$, which is also related to the original size $r$. This part is treated in section \ref{SecEvolMol2}. 
\end{itemize}
To conclude it is enough to iterate the second step as many times as necessary to get rid of the dependence of the times $s_0, s_1,...$ from the molecule's size. This way we obtain the $L^q$ control needed for all time $T_0<t<T$.

\subsection{Small time molecule's evolution: First step}\label{SecEvolMol1}
The following theorem shows how the molecular properties are deformed with the evolution for a small time $s_0$.
\begin{Theoreme}\label{SmallGeneralisacion} Set $\sigma$, $\gamma$ and $\omega$ three real numbers such that $\frac{N}{N+1}<\sigma<1$, $\gamma=N(\frac{1}{\sigma}-1)$ and $0<\gamma<\omega< 1$. Let $\psi(x,s_0)$ be a solution of the problem
\begin{equation}\label{SmallEvolution}
\left\lbrace
\begin{array}{rl}
\partial_{s_0} \psi(x,s_0)=& -\nabla\cdot(v\, \psi)(x,s_0)-\mathcal{J}^{1/2}\psi(x,s_0)\\[5mm]
\psi(x,0)=& \psi_0(x)\\[5mm]
div(v)=&0 \quad \mbox{and }\; v\in L^{\infty}([0,T];bmo(\mathbb{G}))\quad \mbox{with } \underset{s_0\in [0,T]}{\sup}\; \|v(\cdot,s_0)\|_{bmo}\leq \mu
\end{array}
\right.
\end{equation}
If $\psi_0$ is a small $r$-molecule in the sense of definition \ref{DefMolecules} for the local Hardy space $h^\sigma(\mathbb{G})$, then there exists a positive constant $K=K(\mu)$ big enough and a positive constant $\epsilon$ such that for all $0< s_0 \leq\epsilon r$ small we have the following estimates
\begin{eqnarray}
\int_{\mathbb{G}}|\psi(x,s_0)|\|x\cdot x^{-1}(s_0)\|^\omega dx &\leq &(r+Ks_0)^{\omega-\gamma}  \label{SmallConcentration}\\
\|\psi(\cdot, s_0)\|_{L^\infty}&\leq & \frac{1}{\big(r+Ks_0\big)^{N+\gamma}}\label{SmallLinftyevolution}\\
\|\psi(\cdot, s_0)\|_{L^1} &\leq & \frac{v_N}{\big(r+Ks_0\big)^{\gamma}}\label{SmallL1evolution}
\end{eqnarray}
where $v_N$ denotes the volume of the unit ball.\\ 

The new molecule's center $x(s_0)$ used in formula (\ref{SmallConcentration}) is fixed by 
\begin{equation}\label{Defpointx_0}
\left\lbrace
\begin{array}{rl}
x'(s_0)=& \overline{v}_{B_r}=\frac{1}{|B_r|}\displaystyle{\int_{B_r}}v(y,s_0)dy \quad \mbox{where } B_r=B(x(s_0),r).\\[5mm]
x(0)=& x_0.
\end{array}
\right.
\end{equation}  
\end{Theoreme}
\begin{Remarque}
\begin{itemize}
\item[]
\item[1)] The definition of the point $x(s_0)$ given by (\ref{Defpointx_0}) reflects the molecule's center transport using velocity $v$.
\item[2)] With estimates (\ref{SmallLinftyevolution}) and (\ref{SmallL1evolution}) at our disposal we have 
\begin{equation*}
\|\psi(\cdot, s_0)\|_{L^q}\leq C (r+Ks_0)^{-N(1-1/q)-\gamma}.
\end{equation*}
\item[3)] Remark that it is enough to treat the case $0<(r+Ks_0)<1$ since $s_0$ is small: otherwise the $L^q$ control will be trivial for time $s_0$ and beyond: we only need to apply the maximum principle.
\end{itemize}
\end{Remarque}
\textbf{\textit{Proof of the Theorem \ref{SmallGeneralisacion}.}} 
We will follow the next scheme: first we prove the small Concentration condition (\ref{SmallConcentration}) and then we prove the Height condition (\ref{SmallLinftyevolution}). Once we have these two  conditions, the $L^1$ estimate (\ref{SmallL1evolution}) will follow easily. 
\subsubsection*{1) Small time Concentration condition} 
Let us write $\Omega_0(x)=\|x\cdot x^{-1}(s_0)\|^{\omega}$ and $\psi(x)=\psi_+(x)-\psi_-(x)$ where the functions $\psi_{\pm}(x)\geq 0$ have disjoint support. We will note $\psi_\pm(x,s_0)$ two solutions of (\ref{SmallEvolution}) with $\psi_\pm(x,0)=\psi_\pm(x)$. At this point, we use the positivity principle, thus by linearity we have
$$|\psi(x,s_0)|=|\psi_+(x,s_0)-\psi_-(x,s_0)|\leq \psi_+(x,s_0)+\psi_-(x,s_0)$$ and we can write
$$\int_{\mathbb{G}}|\psi(x,s_0)|\Omega_0(x)dx\leq\int_{\mathbb{G}}\psi_+(x,s_0)\Omega_0(x)dx+\int_{\mathbb{G}}\psi_-(x,s_0)\Omega_0(x)dx$$
so we only have to treat one of the integrals on the right side above. We have:
\begin{eqnarray*}
I&=&\left|\partial_{s_0} \int_{\mathbb{G}}\Omega_0(x)\psi_+(x,s_0)dx\right|\\
&=&\left|\int_{\mathbb{G}}\partial_{s_0} \Omega_0(x)\psi_+(x,s_0)+\Omega_0(x)\left[-\nabla\cdot(v\, \psi_+(x,s_0))-\mathcal{J}^{1/2}\psi_+(x,s_0)\right]dx\right|\\
&=&\left|\int_{\mathbb{G}}-\nabla\Omega_0(x)\cdot x'(s_0)\psi_+(x,s_0)+\Omega_0(x)\left[-\nabla\cdot(v\, \psi_+(x,s_0))-\mathcal{J}^{1/2}\psi_+(x,s_0)\right]dx\right|\\
&
\end{eqnarray*}
Using the fact that $v$ is divergence free, we obtain
\begin{equation*}
I=\left|\int_{\mathbb{G}}\nabla\Omega_0(x)\cdot(v-x'(s_0))\psi_+(x,s_0)-\Omega_0(x)\mathcal{J}^{1/2}\psi_+(x,s_0)dx\right|.
\end{equation*}
Since the operator $\mathcal{J}^{1/2}$ is symmetric and using the definition of $x'(s_0)$ given in (\ref{Defpointx_0}) we have
\begin{equation}\label{smallEstrella}
I\leq c \underbrace{\int_{\mathbb{G}}\|x\cdot x^{-1}(s_0)\|^{\omega-1}|v-\overline{v}_{B_r}| |\psi_+(x,s_0)|dx}_{I_1} + c\underbrace{\int_{\mathbb{G}}\big|\mathcal{J}^{1/2}\Omega_0(x)\big|\, |\psi_+(x,s_0)|dx}_{I_2}.
\end{equation}
We will study separately each of the integrals $I_1$ and $I_2$:
\begin{Lemme}\label{Lemme1} For integral $I_1$ above we have the estimate $I_1\leq C \mu \; r^{\omega-1-\gamma}$.
\end{Lemme}
\begin{Lemme}\label{Lemme2} For integral $I_2$ in inequality (\ref{smallEstrella}) we have the inequality $I_2\leq C r^{\omega-1-\gamma}$.
\end{Lemme}
Using these lemmas and getting back to estimate (\ref{smallEstrella}) we have
$$\left|\partial_{s_0} \int_{\mathbb{G}}\Omega_0(x)\psi_+(x,s_0)dx\right| \leq  C(\mu+1)\; r^{\omega-1-\gamma}$$
This last estimation is compatible with the estimate (\ref{SmallConcentration}) for $0\leq s_0\leq \epsilon r$ small enough: just fix $K$ such that
\begin{equation}\label{SmallConstants}
C\left(\mu+1\right)\leq K(\omega-\gamma).
\end{equation}
Indeed, since the time $s_0$ is very small, we can linearize the formula $(r+Ks_0)^{\omega-\gamma}$ in the right-hand side of (\ref{SmallConcentration}) in order to obtain
\begin{equation*}
\phi=(r+Ks_0)^{\omega-\gamma} \thickapprox r^{\omega-\gamma}\left(1+[K(\omega-\gamma)]\frac{s_0}{r}\right).
\end{equation*}
Finally, taking the derivative with respect to $s_0$ in the above expression we have $\phi' \thickapprox r^{\omega-1-\gamma}K(\omega-\gamma)$ and with condition (\ref{SmallConstants}), the small time Concentration condition (\ref{SmallConcentration}) follows.\\

We prove now the Lemmas \ref{Lemme1} and \ref{Lemme2}; but before, we will need the following result
\begin{Lemme}\label{PropoBMO1} Let $f\in bmo(\mathbb{G})$, then
\begin{enumerate}
\item[1)] for all $1<p<+\infty$, $f$ is locally in $L^p$ and $\frac{1}{|B|}\displaystyle{\int_{B}|f(x)-f_B|^pdx}\leq C \|f\|_{bmo}^p$
\item[2)] for all $k\in \mathbb{N}$, we have $|f_{2^k B}-f_B|\leq Ck \|f\|_{bmo}$ where $2^kB=B(x,2^k R)$ is a ball centered at a point $x$ of radius $2^kR$. 
\end{enumerate}
\end{Lemme}
For a proof of these results see \cite{Stein2}.\\

\textit{\textbf{Proof of the Lemma \ref{Lemme1}.}} We begin by considering the space $\mathbb{G}$ as the union of a ball with dyadic coronas centered around $x(s_0)$, more precisely we set $\mathbb{G}=B_r\cup \bigcup_{k\geq 1}E_k$ where
\begin{equation}\label{SmallDecoupage}
B_r= \{x\in \mathbb{G}: \|x\cdot x^{-1}(s_0)\|\leq r\} \quad \mbox{and}\quad E_k= \{x\in \mathbb{G}: r2^{k-1}<\|x\cdot x^{-1}(s_0)\|\leq r 2^{k}\} \quad \mbox{for } k\geq 1,
\end{equation}
\begin{enumerate}
\item[(i)] \underline{Estimations over the ball $B_r$}. Applying Hölder's inequality to the integral $I_{1,B_r}$ we obtain
\begin{eqnarray}
I_{1,B_r}=\int_{B_r}\|x\cdot x^{-1}(s_0)\|^{\omega-1}|v-\overline{v}_{B_r}| |\psi_+(x,s_0)|dx & \leq &\underbrace{\| \|x\cdot x^{-1}(s_0)\|^{\omega-1}\|_{L^p(B_r)}}_{(1)} \label{Equa1} \\
& \times &\underbrace{\|v-\overline{v}_{B_r}\|_{L^z(B_r)}}_{(2)}\underbrace{\|\psi_+(\cdot, s_0)\|_{L^q(B_r)}}_{(3)}\nonumber
\end{eqnarray}
where $\frac{1}{p}+\frac{1}{z}+\frac{1}{q}=1$ and $p,z,q> 1$. We treat each of the previous terms separately:
\begin{enumerate}
\item[$\bullet$] First observe that for $1<p<N/(1-\omega)$ we have for the term $(1)$ above:
\begin{equation*}
\|\|x\cdot x^{-1}(s_0)\|^{\omega-1}\|_{L^p(B_r)}\leq C r^{N/p+\omega-1}.
\end{equation*}
\item[$\bullet$] By hypothesis $v(\cdot, s_0)\in bmo$ and applying Lemma \ref{PropoBMO1} we have $\|v-\overline{v}_{B_r}\|_{L^z(B_r)}\leq C|B_r|^{1/z}\|v(\cdot,s_0)\|_{bmo}$. Since $\underset{s_0\in [0,T]}{\sup}\; \|v(\cdot,s_0)\|_{bmo}\leq \mu$, we write for the term $(2)$:
\begin{equation*}
\|v-\overline{v}_{B_r}\|_{L^z(B_r)}\leq C\mu\; r^{N/z}.
\end{equation*}
\item[$\bullet$] Finally for $(3)$ by the maximum principle we have $\|\psi_+(\cdot, s_0)\|_{L^q(B_r)}\leq  \|\psi_+(\cdot, 0)\|_{L^q}$; hence using the fact that $\psi_0$ is an $r$-molecule and remark \ref{Remark2}-2) we obtain
\begin{equation*}
\|\psi_+(\cdot, s_0)\|_{L^q(B_r)}\leq  C r^{-N(1-1/q)-\gamma}.
\end{equation*}
\end{enumerate}
We combine all these inequalities together in order to obtain the following estimation for (\ref{Equa1}):
\begin{equation}\label{Bola1}
I_{1,B_r}\leq C\mu\; r^{\omega-1-\gamma}.
\end{equation}
\item[(ii)] \underline{Estimations for the dyadic corona $E_k$}. Let us note $I_{1,E_k}$ the integral 
$$I_{1,E_k}=\int_{E_k}\|x\cdot x^{-1}(s_0)\|^{\omega-1}|v-\overline{v}_{B_r}| |\psi_+(x,s_0)|dx.$$
Since over $E_k$ we have\footnote{recall that $0<\gamma<\omega<1$.} $\|x\cdot x^{-1}(s_0)\|^{\omega-1}\leq C 2^{k(\omega-1)}r^{\omega-1}$ we write
\begin{eqnarray*}
I_{1,E_k}&\leq & C2^{k(\omega-1)}r^{\omega-1}\left(\int_{E_k}|v-\overline{v}_{B_{r2^k}}| |\psi_+(x,s_0)|dx+\int_{E_k}|\overline{v}_{B_r}-\overline{v}_{B_{r2^k}}| |\psi_+(x,s_0)|dx\right)
\end{eqnarray*}
where we noted $B_{r2^k}=B(x(s_0),r2^k)$, then
\begin{eqnarray*}
I_{1,E_k} &\leq& C2^{k(\omega-1)}r^{\omega-1}\left(\int_{B_{r2^k}}|v-\overline{v}_{B_{r2^k}}| |\psi_+(x,s_0)|dx+\int_{B_{r2^k}}|\overline{v}_{B_r}-\overline{v}_{B_{r2^k}}| |\psi_+(x,s_0)|dx\right).
\end{eqnarray*}
Now, since $v(\cdot, s_0)\in bmo(\mathbb{G})$, using the Lemma \ref{PropoBMO1} we have $|\overline{v}_{B_r}-\overline{v}_{B_{r2^k}}| \leq Ck\|v(\cdot,s_0)\|_{bmo}\leq Ck\mu$ and we can write
\begin{eqnarray*}
I_{1,E_k}&\leq &C2^{k(\omega-1)}r^{\omega-1}\left(\int_{B_{r2^k}}|v-\overline{v}_{B_{r2^k}}| |\psi_+(x,s_0)|dx+Ck\mu\|\psi_+(\cdot,s_0)\|_{L^1}\right)\\[5mm]
&\leq & C2^{k(\omega-1)}r^{\omega-1}\left(\|\psi_+(\cdot,s_0)\|_{L^{a_0}} \|v-\overline{v}_{B_{r2^k}}\|_{L^{\frac{a_0}{a_0-1}}} +Ck\mu\;  r^{-\gamma}\right)
\end{eqnarray*}
where we used Hölder's inequality with $1<a_0<\frac{N}{N+(\omega-1)}$ and maximum principle for the last term above. Using again the properties of $bmo$ spaces we have
$$I_{1,E_k}\leq  C2^{k(\omega-1)}r^{\omega-1}\left(\|\psi_+(\cdot,0)\|_{L^1}^{1/a_0}\|\psi_+(\cdot,0)\|_{L^\infty}^{1-1/a_0} |B_{r2^k}|^{1-1/a_0}\|v(\cdot,s)\|_{bmo} +Ck\mu r^{-\gamma}\right).$$
Let us now apply the estimates given by hypothesis for $\|\psi_+(\cdot,0)\|_{L^1}$, $\|\psi_+(\cdot,0)\|_{L^\infty}$ and $\|v(\cdot,s_0)\|_{bmo}$  to obtain
$$I_{1,E_k}\leq  C2^{k(N-N/a_0+\omega-1)}r^{\omega-1-\gamma}\mu +C2^{k(\omega-1)}k\mu\;  r^{\omega-1-\gamma}.$$
Since $1<a_0<\frac{N}{N+(\omega-1)}$, we have $N-N/a_0+(\omega-1)<0$, so that, summing over each dyadic corona $E_k$, we have
\begin{equation}\label{Coronak}
\sum_{k\geq 1}I_{1,E_k}\leq C\mu\; r^{\omega-1-\gamma}.
\end{equation}
\end{enumerate}
Finally, gathering together the estimations (\ref{Bola1}) and (\ref{Coronak}) we obtain the desired conclusion.\hfill$\blacksquare$\\

\textbf{\textit{Proof of the Lemma \ref{Lemme2}.}} As for the Lemma \ref{Lemme1}, we consider $\mathbb{G}$ as the union of a ball with dyadic coronas centered on $x(s_0)$ (cf. (\ref{SmallDecoupage})). 
\begin{enumerate}
\item[(i)] \underline{Estimations over the ball $B_r$}. We write, using the maximum principle and the hypothesis for $\|\psi_+(\cdot, 0)\|_{L^\infty}$:
\begin{eqnarray*}
I_{2,B_r}&=&\int_{B_r}\big|\mathcal{J}^{1/2}(\|x\cdot x^{-1}(s_0)\|^{\omega})\big||\psi_+(x,s_0)|dx \leq \|\psi_+(\cdot, s_0)\|_{L^{\infty}}\int_{B_r}\|x\cdot x^{-1}(s_0)\|^{\omega-1}dx\nonumber\\
&\leq & C r^{-N-\gamma} r^{N+\omega-1}=Cr^{\omega-1-\gamma}.
\end{eqnarray*}
\item[(ii)]  \underline{Estimations for the dyadic corona $E_k$}. 
\begin{eqnarray*}
I_{2, E_k}&=&\int_{E_k}|\mathcal{J}^{1/2}(\|x\cdot x^{-1}(s_0)\|^{\omega})| |\psi_+(x,s_0)|dx  \leq \left(\underset{x\in E_k}{\sup}\|x\cdot x^{-1}(s_0)\|^{\omega-1}\right) \|\psi_+(\cdot, s_0)\|_{L^1}\\
&\leq &C r^{-\gamma} \big(2^k r\big)^{\omega-1}=C r^{\omega-1-\gamma} 2^{-k(1-\omega)}
\end{eqnarray*}
Since $0<\gamma<\omega<1$, summing over $k\geq 1$, we obtain 
\begin{equation*}
\sum_{k\geq 1}I_{2,E_k}\leq C r^{\omega-1-\gamma}.
\end{equation*}
\end{enumerate}
In order to finish the proof of Lemma \ref{Lemme2} we combine together the estimates (i) and (ii).\hfill$\blacksquare$

\subsubsection*{2) Small time Height condition}
We treat now the Height condition (\ref{SmallLinftyevolution}) and for this we will give a sligthly different proof of the maximum principle of A. C\'ordoba \& D. C\'ordoba. Indeed, the following proof only relies on the Concentration condition. \\

Assume that molecules we are working with are smooth enough. Following an idea of \cite{Cordoba} (section 4 p.522-523), we will note $\overline{x}$ the point of $\mathbb{G}$ such that $\psi(\overline{x},s_0)=\|\psi(\cdot,s_0)\|_{L^\infty}$. Thus we can write:
\begin{equation}\label{Infty1}
\frac{d}{ds_0}\|\psi(\cdot,s_0)\|_{L^\infty}\leq  -\int_{\{\|\overline{x}\cdot y^{-1}\|<1\}}\frac{\psi(\overline{x},s_0)-\psi(y,s_0)}{\|\overline{x}\cdot y^{-1}\|^{N+1}}dy\leq 0.
\end{equation}
Let us consider the corona centered in $\overline{x}$ defined by
$$\mathcal{C}(R_1,R_2)=\{y\in \mathbb{G}:R_1\leq\|\overline{x}\cdot y^{-1}\|\leq R_2\}$$ 
where $1>R_2=\rho R_1$ with $\rho >2$ and where $R_1$ will be fixed later. Then:
\begin{equation*}
\int_{\{\|\overline{x}\cdot y^{-1}\|<1\}}\frac{\psi(\overline{x},s_0)-\psi(y,s_0)}{\|\overline{x}\cdot y^{-1}\|^{N+1}}dy\geq \int_{\mathcal{C}(R_1,R_2)}\frac{\psi(\overline{x},s_0)-\psi(y,s_0)}{\|\overline{x}\cdot y^{-1}\|^{N+1}}dy.
\end{equation*}
Define the sets $B_1$ and $B_2$ by $B_1=\{y\in \mathcal{C}(R_1,R_2): \psi(\overline{x},s_0)-\psi(y,s_0)\geq \frac{1}{2}\psi(\overline{x},s_0)\}$ and $B_2=\{y\in \mathcal{C}(R_1,R_2): \psi(\overline{x},s_0)-\psi(y,s_0)< \frac{1}{2}\psi(\overline{x},s_0)\}$ such that $\mathcal{C}(R_1,R_2)=B_1\cup B_2$. \\

We obtain the inequalities
\begin{eqnarray*}
\int_{\mathcal{C}(R_1,R_2)}\frac{\psi(\overline{x},s_0)-\psi(y,s_0)}{ \|\overline{x}\cdot y^{-1} \|^{N+1}}dy &\geq &\int_{B_1}\frac{\psi(\overline{x},s_0)-\psi(y,s_0)}{\|\overline{x}\cdot y^{-1}\|^{N+1}}dy \geq \frac{\psi(\overline{x},s_0)}{2R_2^{N+1}}|B_1|=\frac{\psi(\overline{x},s_0)}{2R_2^{N+1}}\left(|\mathcal{C}(R_1,R_2)|-|B_2|\right).
\end{eqnarray*}
Since $R_2=\rho R_1$ one has
\begin{equation}\label{Infty3}
\int_{\mathcal{C}(R_1,R_2)}\frac{\psi(\overline{x},s_0)-\psi(y,s_0)}{\|\overline{x}\cdot y^{-1}\|^{N+1}}dy\geq \frac{\psi(\overline{x},s_0)}{2\rho^{N+1}R_1^{N+1}}\bigg(v_N(\rho^N -1)R_1^N-|B_2|\bigg)
\end{equation}
where $v_N$ denotes the volume of the unit ball.\\ 

To continue, we need to estimate the quantity $|B_2|$ in the right-hand side of (\ref{Infty3}) in terms of $\psi(\overline{x},s_0)$ and $R_1$. We will distinguish two cases:
\begin{enumerate}
\item[1)] if $\|\overline{x}\cdot x^{-1}(s_0)\|>2R_2$ or $\|\overline{x}\cdot x^{-1}(s_0)\|<R_1/2$ then
\begin{equation}\label{FormEstima1}
C_1(r+Ks_0)^{\omega-\gamma}\psi(\overline{x},s_0)^{-1}R_1^{-\omega}\geq |B_2|
\end{equation}
\item[2)] if $R_1/2\leq \|\overline{x}\cdot x^{-1}(s_0)\|\leq 2R_2$ then 
\begin{equation}\label{FormEstima2}
\big(C_2 (r+Ks_0)^{\omega-\gamma}R_1^{n-\omega}\psi(\overline{x},s_0)^{-1}\big)^{1/2}\geq |B_2|.
\end{equation}
\end{enumerate}
For these two estimates, our starting point is the Concentration condition :
\begin{eqnarray}
(r+Ks_0)^{\omega-\gamma}&\geq &\int_{\mathbb{G}}|\psi(y,s_0)|\|y\cdot x^{-1}(s_0)\|^{\omega}dy\nonumber\\ 
&\geq & \int_{B_2}|\psi(y,s_0)|\|y\cdot x^{-1}(s_0)\|^{\omega}dy \geq \frac{\psi(\overline{x},s_0)}{2}\int_{B_2}\|y\cdot x^{-1}(s_0)\|^{\omega}dy.\label{EstimationB2}
\end{eqnarray}
We just need to estimate the last integral following the cases given above. Indeed, if $\|\overline{x}\cdot x^{-1}(s_0)\|>2R_2$  then we have
$$\underset{y\in B_2\subset \mathcal{C}(R_1,R_2)}{\min}\|y\cdot x^{-1}(s_0)\|^{\omega}\geq R_2^{\omega}=\rho^{\omega}R_1^{\omega}$$
while if $\|\overline{x}\cdot x^{-1}(s_0)\|<R_1/2$, one has
$$\underset{y\in B_2\subset \mathcal{C}(R_1,R_2)}{\min}\|y\cdot x^{-1}(s_0)\|^{\omega}\geq \frac{R_1^{\omega}}{2^\omega}.$$
Applying these results to (\ref{EstimationB2}) we obtain $(r+Ks_0)^{\omega-\gamma}\geq \frac{\psi(\overline{x},s_0)}{2} \rho^{\omega} R_1^{\omega}|B_2|$ and  $(r+Ks_0)^{\omega-\gamma}\geq \frac{\psi(\overline{x},s_0)}{2} \frac{R_1^{\omega}}{2^\omega}|B_2|$, and since $\rho>2$ we have the first desired estimate
\begin{equation*}
\frac{C_1 (r+Ks_0)^{\omega-\gamma}}{\psi(\overline{x},s_0) R_1^{\omega}} \geq \frac{2 (r+Ks_0)^{\omega-\gamma}}{\rho^{\omega}\psi(\overline{x},s_0) R_1^{\omega}} \geq |B_2| \qquad \mbox{with } C_1=2^{1+\omega}.
\end{equation*}
For the second case, since $R_1/2\leq \|\overline{x}\cdot x^{-1}(s_0)\|\leq 2R_2$, we can write using the Cauchy-Schwarz inequality
\begin{equation}\label{HolderInver}
\int_{B_2}\|y\cdot x^{-1}(s_0)\|^{\omega}dy\geq |B_2|^2\left(\int_{B_2}\|y\cdot x^{-1}(s_0)\|^{-\omega}dy\right)^{-1}
\end{equation}
Now, observe that in this case we have $B_2\subset B(x(s_0), 5R_2)$ and then
$$\int_{B_2}\|y\cdot x^{-1}(s_0)\|^{-\omega}dy\leq \int_{B(x(s_0), 5 R_2)}\|y\cdot x^{-1}(s_0)\|^{-\omega}dy\leq v_N (5\rho R_1)^{N-\omega}.$$
Getting back to(\ref{HolderInver}) we have
$$\int_{B_2}\|y\cdot x^{-1}(s_0)\|^{\omega}dy\geq |B_2|^2 v_N^{-1} (5 \rho R_1)^{-N+\omega}$$
and we use this estimate in (\ref{EstimationB2}) to obtain
\begin{equation*}
\frac{C_2(r+Ks_0)^{\omega/2-\gamma/2} R_1^{N/2-\omega/2}}{\psi(\overline{x},s_0)^{1/2}}\geq |B_2|, \quad \mbox{where } C_2=(2\times 5^{N-\omega} v_N\rho^{N-\omega})^{1/2}.\\[5mm]
\end{equation*}
Now, with estimates (\ref{FormEstima1}) and (\ref{FormEstima2}) at our disposal we can write
\begin{enumerate}
\item[(i)] if $\|\overline{x}\cdot x^{-1}(s_0)\|>2R_2$ or $\|\overline{x}\cdot x^{-1}(s_0)\|<R_1/2$ then
\begin{equation*}
\int_{\mathcal{C}(R_1,R_2)}\frac{\psi(\overline{x},s_0)-\psi(y,s_0)}{\|\overline{x}\cdot y^{-1}\|^{N+1}}dy\geq  \frac{\psi(\overline{x},s_0)}{2\rho^{N+1}R_1^{N+1}}\bigg(v_N(\rho^N -1)R_1^n-\frac{C_1(r+Ks_0)^{\omega-\gamma}}{\psi(\overline{x},s_0)} R_1^{-\omega}\bigg)
\end{equation*}
\item[(ii)] if $R_1/2\leq \|\overline{x}\cdot x^{-1}(s_0)\|\leq 2R_2$
\begin{equation*}
\int_{\mathcal{C}(R_1,R_2)}\frac{\psi(\overline{x},s_0)-\psi(y,s_0)}{\|\overline{x}\cdot y^{-1}\|^{N+1}}dy\geq  \frac{\psi(\overline{x},s_0)}{2\rho^{N+1}R_1^{N+1}}\bigg(v_N(\rho^N -1)R_1^N-\frac{C_2 (r+Ks_0)^{\omega/2-\gamma/2}R_1^{N/2-\omega/2}}{\psi(\overline{x},s_0)^{1/2}}\bigg)
\end{equation*}
\end{enumerate}
If we set $R_1=(r+Ks_0)^{\frac{(\omega-\gamma)}{N+\omega}}\psi(\overline{x},s_0)^{\frac{-1}{N+\omega}}$ and if $\rho$ is big enough such that the expressions in brackets above are positive, we obtain for cases (i) and (ii) the following estimate for (\ref{Infty3}):
\begin{equation*}
\int_{\mathcal{C}(R_1,R_2)}\frac{\psi(\overline{x},s_0)-\psi(y,s_0)}{\|\overline{x}\cdot y^{-1}\|^{N+1}}dy\geq C (r+Ks_0)^{-\frac{(\omega-\gamma)}{N+\omega}} \psi(\overline{x},s_0)^{1+\frac{1}{N+\omega}}
\end{equation*}
where $C=C(N,\rho)=\frac{v_N (\rho^N-1)-\sqrt{2v_N}(5\rho)^{\frac{N-\omega}{2}}}{2\rho^{N+1}}<1$ is a small positive constant. Now, and for all possible cases considered before, we have the following estimate for (\ref{Infty1}):
$$\frac{d}{ds_0}\|\psi(\cdot,s_0)\|_{L^\infty}\leq -C(r+Ks_0)^{-\frac{(\omega-\gamma)}{N+\omega}} \|\psi(\cdot,s_0)\|_{L^\infty}^{1+\frac{1}{N+\omega}}.$$
In order to solve this problem, it is enough to remark that if $\|\psi(\cdot,s_0)\|_{L^\infty}\leq(r+Ks_0)^{-(N+\gamma)}$, then $\|\psi(\cdot,s_0)\|_{L^\infty}$ satisfies the previous inequality. Indeed, we have
\begin{eqnarray*}
\frac{d}{ds_0}\|\psi(\cdot,s_0)\|_{L^\infty}&\leq &-K(N+\gamma)(r+Ks_0)^{-(N+\gamma)-1}\\
&\leq & -C(r+Ks_0)^{-(N+\gamma)-1}=-C(r+Ks_0)^{-\frac{(\omega-\gamma)}{N+\omega}}(r+Ks_0)^{-(N+\gamma)(1+\frac{1}{N+\omega})}\\
&\leq & -C(r+Ks_0)^{-\frac{(\omega-\gamma)}{N+\omega}}\|\psi(\cdot,s_0)\|_{L^\infty}^{1+\frac{1}{N+\omega}}
\end{eqnarray*}
Furthermore with the initial data $\|\psi(\cdot, 0)\|_{L^{\infty}}\leq r^{-N-\gamma}$, we obtain that this solution is unique.
\subsubsection*{3) Small time $L^q$ estimate}
This last condition is an easy consequence of the previous computations. Indeed: we write
\begin{eqnarray}
\int_{\mathbb{G}}|\psi(x,s_0)|dx&=&\int_{\{\|x\cdot x^{-1}(s_0)\|< D\}}|\psi(x,s_0)|dx+\int_{\{\|x\cdot x^{-1}(s_0)\|\geq D\}}|\psi(x,s_0)|dx\nonumber\\
&\leq & v_N D^N \|\psi(\cdot, s_0)\|_{L^\infty}+D^{-\omega}\int_{\mathbb{G}}|\psi(x,s_0)|\|x\cdot x^{-1}(s_0)\|^\omega dx\nonumber
\end{eqnarray}
Now using the Concentration condition and the Height condition one has:
\begin{eqnarray*}
\int_{\mathbb{G}}|\psi(x,s_0)|dx&\leq & v_N \frac{D^N }{\left(r+Ks_0\right)^{N+\omega}}  +D^{-\omega}(r+Ks_0)^{\omega-\gamma}
\end{eqnarray*}
To continue, it is enough to choose correctly the real parameter $D$ to obtain
$$\|\psi(\cdot, s_0)\|_{L^1}\leq \frac{v_N}{\big(r+Ks_0\big)^{\gamma}}.$$
Once we have the $L^\infty$ and the $L^1$ bounds, the estimate for the norm $L^q$ is inmediate and Theorem \ref{SmallGeneralisacion} is now completely proven. \hfill$\blacksquare$
\subsection{Molecule's evolution: Second step}\label{SecEvolMol2}
In the previous section we have obtained deformed molecules after a very small time $s_0$. The next theorem shows us how to obtain similar profiles in the inputs and the outputs in order to perform an iteration in time. 
\begin{Theoreme}\label{Generalisacion} Set $\gamma$ and $\omega$ two real numbers such that $0<\gamma<\omega<1$. Let $0< s_1\leq T$ and let $\psi(x,s_1)$ be a solution of the problem
\begin{equation}\label{Evolution}
\left\lbrace
\begin{array}{rl}
\partial_{s_1} \psi(x,s_1)=& -\nabla\cdot(v\, \psi)(x,s_1)-\mathcal{J}^{1/2}\psi(x,s_1)\\[5mm]
\psi(x,0)=& \psi(x,s_0)  \qquad \qquad \mbox{with } s_0>0\\[5mm]
div(v)=&0 \quad \mbox{and }\; v\in L^{\infty}([0,T];bmo(\mathbb{G}))\quad \mbox{with } \underset{s_1\in [s_0,T]}{\sup}\; \|v(\cdot,s_1)\|_{bmo}\leq \mu
\end{array}
\right.
\end{equation}
If $\psi(x,s_0)$ satisfies the three following conditions
\begin{eqnarray*}
\int_{\mathbb{G}}|\psi(x,s_0)|\|x \cdot x^{-1}(s_0)\|^\omega dx \leq (r+Ks_0)^{\omega-\gamma}; \quad \|\psi(\cdot, s_0)\|_{L^\infty}\leq \frac{1}{\left(r+ Ks_0\right)^{N+\gamma }}; \quad \|\psi(\cdot, s_0)\|_{L^1} \leq  \frac{v_N}{\big(r+Ks_0\big)^{\gamma }}
\end{eqnarray*}
where $K=K(\mu)$ is given by (\ref{SmallConstants}) and $s_0$ is such that $(r+Ks_0)<1$. Then for all $0< s_1\leq\epsilon r$ small, we have the following estimates
\begin{eqnarray}
\int_{\mathbb{G}}|\psi(x,s_1)|\|x \cdot x^{-1}(s_1)\|^\omega dx &\leq &(r+K(s_0+s_1))^{\omega-\gamma}  \label{Concentration2}\\
\|\psi(\cdot,s_1)\|_{L^\infty}&\leq & \frac{1}{\left(r+K(s_0+s_1)\right)^{N+\gamma }}\label{Linftyevolutionnotsmalltime}\\
\|\psi(\cdot,s_1)\|_{L^1} &\leq & \frac{v_N}{\big(r+K(s_0+s_1)\big)^{\gamma }} \label{L1evolutionsmalltime}
\end{eqnarray}
\end{Theoreme}

\begin{Remarque}
\begin{itemize}
\item[]
\item[1)] Since $s_1$ is small and $(r+Ks_0)<1$, we can without loss of generality assume that $(r+K(s_0+s_1))<1$: otherwise, by the maximum principle there is nothing to prove.
\item[2)] The new molecule's center $x(s_1)$ used in formula (\ref{Concentration2}) is fixed by 
\begin{equation}\label{Defpointx_s}
\left\lbrace
\begin{array}{rl}
x'(s_1)=& \overline{v}_{B_{f_1}}=\frac{1}{|B_{f_1}|}\displaystyle{\int_{B_{f_1}}}v(y,s_1)dy\\[5mm]
x(0)=& x(s_0).
\end{array}
\right.
\end{equation}  
And here we noted $B_{f_1}=B(x(s_1),f_1)$ with $f_1$ a real valued function given by 
\begin{equation}\label{DefiFunctionF}
f_1=(r+Ks_0).
\end{equation}
Note that by remark $1)$ above we have $0<f_1<1$.
\item[3)] We recall that the wished $L^q$ bound is given by interpolating (\ref{Linftyevolutionnotsmalltime}) and (\ref{L1evolutionsmalltime}).
\end{itemize}
\end{Remarque}
\textit{\textbf{Proof of the Theorem  \ref{Generalisacion}.}}
We will follow the same scheme as before:  we first prove the Concentration condition (\ref{Concentration2}), with this estimate at hand we will control the $L^\infty$ decay and then we will obtain the suitable $L^1$ control.
\subsubsection*{1) The Concentration condition}
The calculations are very similar of those performed before: the only difference stems from the initial data and the definition of the center $x(s_1)$. So, let us write $\Omega_1(x)=\|x\cdot x^{-1}(s_1)\|^{\omega}$ and $\psi(x)=\psi_+(x)-\psi_-(x)$ where the functions $\psi_{\pm}(x)\geq 0$ have disjoint support. Thus, by linearity and using the positivity principle we have  
$$|\psi(x,s_1)|=|\psi_+(x,s_1)-\psi_-(x,s_1)|\leq \psi_+(x,s_1)+\psi_-(x,s_1)$$ and we can write
$$\int_{\mathbb{G}}|\psi(x,s_1)\Omega_1(x)dx\leq\int_{\mathbb{G}}\psi_+(x,s_1)\Omega_1(x)dx+\int_{\mathbb{G}}\psi_-(x,s_1)\Omega_1(x)dx$$
so we only have to treat one of the integrals on the right-hand side above. We have:
\begin{eqnarray*}
I&=&\left|\partial_{s_1} \int_{\mathbb{G}}\Omega_1(x)\psi_+(x,s_1)dx\right|=\left|\int_{\mathbb{G}}-\nabla\Omega_1(x)\cdot x'(s_1)\psi_+(x,s_1)+\Omega_1(x)\left[-\nabla\cdot(v\, \psi_+(x,s_1))-\mathcal{J}^{1/2}\psi_+(x,s_1)\right]dx\right|
\end{eqnarray*}
Using the fact that $v$ is divergence free, we obtain
\begin{equation*}
I=\left|\int_{\mathbb{G}}\nabla\Omega_1(x)\cdot(v-x'(s_1))\psi_+(x,s_1)-\Omega_1(x)\mathcal{J}^{1/2}\psi_+(x,s_1)dx\right|.
\end{equation*}
Finally, using the definition of $x'(s_1)$ given in (\ref{Defpointx_s}) and replacing $\Omega_1(x)$ by $\|x\cdot x^{-1}(s_1)\|^{\omega}$ in the first integral we obtain
\begin{equation}\label{Estrella}
I\leq c \underbrace{\int_{\mathbb{G}}\|x \cdot x^{-1}(s_1)\|^{\omega-1}|v-\overline{v}_{B_{f_1}}| |\psi_+(x,s_1)|dx}_{I_1} + c\underbrace{\int_{\mathbb{G}}\|x \cdot x^{-1}(s_1)\|^{\omega-1}|\psi_+(x,s_1)|dx}_{I_2}.
\end{equation}
We will study separately each of the integrals $I_1$ and $I_2$ in the next lemmas:
\begin{Lemme}\label{Lemme41} For integral $I_1$ we have the estimate $I_1\leq C \mu\big(r+Ks_0\big)^{\omega-\gamma-1}$.
\end{Lemme}
\begin{Lemme}\label{Lemme42}
For integral $I_2$ in the inequality (\ref{Estrella}) we have the following estimate $I_2\leq C\big(r+Ks_0\big)^{\omega-\gamma-1}$.
\end{Lemme}
Using these lemmas and getting back to the estimate (\ref{Estrella}) we have
\begin{equation}\label{FinalEstimate}
\left|\partial_{s_1} \int_{\mathbb{G}}\Omega_1(x)\psi_+(x,s_1)dx\right| \leq  C \left(\mu+1\right)\big(r+Ks_0\big)^{\omega-\gamma-1}
\end{equation}
This estimation is compatible with the estimate (\ref{Concentration2}) for $0\leq s_1\leq \epsilon r$ small enough. Indeed, we can write
$\phi=(r+K(s_0+s_1))^{\omega-\gamma}$ and we linearize this expression with respect to $s_1$:
$$\phi \thickapprox(r+s_0)^{\omega-\gamma}\left(1+K(\omega-\gamma)\frac{s_1}{(r+s_0)}\right)$$
Taking the derivative of $\phi$ with respect to $s_1$ we have $\phi' \thickapprox K(\omega-\gamma) \big(r+Ks_0\big)^{\omega-\gamma-1}$ and with the condition (\ref{SmallConstants}) on $K(\omega-\gamma)$ we obtain that (\ref{FinalEstimate}) is bounded by $\phi'$ and the Concentration condition follows. \\

\textit{\textbf{Proof of the Lemma \ref{Lemme41}}.} We begin by considering the space $\mathbb{G}$ as the union of a ball with dyadic coronas centered on $x(s_1)$, more precisely we set $\mathbb{G}=B_{f_1}\cup\bigcup_{k\geq 1}E_k$ where
\begin{eqnarray}\label{Decoupage}
B_{f_1}&=& \{x\in \mathbb{G}: \|x\cdot x^{-1}(s_1)\|\leq f_1\},\\[5mm]
E_k&=& \{x\in \mathbb{G}: f_12^{k-1}<\|x\cdot x^{-1}(s_1)\|\leq f_1 2^{k}\} \quad \mbox{for } k\geq 1.\nonumber
\end{eqnarray}
\begin{enumerate}
\item[(i)] \underline{Estimations over the ball $B_{f_1}$}. Applying Hölder's inequality on integral $I_1$ we obtain
\begin{eqnarray*}
I_{1,B_{f_1}}=\int_{B_{f_1}}\|x\cdot x^{-1}(s_1)\|^{\omega-1}|v-\overline{v}_{B_{f_1}}| |\psi_+(x,s_1)|dx & \leq &\underbrace{\|\|x\cdot x^{-1}(s_1)\|^{\omega-1}\|_{L^p(B_{f_1})}}_{(1)} \\
& \times &\underbrace{\|v-\overline{v}_{B_{f_1}}\|_{L^z(B_{f_1})}}_{(2)}\underbrace{\|\psi_+(\cdot, s_1)\|_{L^q(B_{f_1})}}_{(3)}
\end{eqnarray*}
where $\frac{1}{p}+\frac{1}{z}+\frac{1}{q}=1$ and $p,z,q> 1$. We treat each of the previous terms separately:
\begin{enumerate}
\item[$\bullet$] Observe that for $1<p<N/(1-\omega)$ we have
\begin{equation*}
\|\|x\cdot x^{-1}(s_1)\|^{\omega-1}\|_{L^p(B_{f_1})}\leq C f_1^{n/p+\omega-1}.
\end{equation*}
\item[$\bullet$] We have $v(\cdot, s_1)\in bmo(\mathbb{G})$, thus $\|v-\overline{v}_{B_{f_1}}\|_{L^z(B_{f_1})}\leq C|B_{f_1}|^{1/z}\|v(\cdot,s_1)\|_{bmo}$. 
Since $\underset{s_1\in [s_0,T]}{\sup}\; \|v(\cdot,s_1)\|_{bmo}\leq \mu$ we write
\begin{equation*}
\|v-\overline{v}_{B_{f_1}}\|_{L^z(B_{f_1})}\leq C f_1^{N/z}\mu.
\end{equation*}
\item[$\bullet$] Finally, by the maximum principle for $L^q$ norms we have $\|\psi_+(\cdot, s_1)\|_{L^q(B_{f_1})}\leq  \|\psi(\cdot, s_0)\|_{L^q}$; hence we obtain
\begin{equation*}
\|\psi_+(\cdot, s_1)\|_{L^q(B_{f_1})}\leq  \|\psi(\cdot, s_0)\|_{L^1}^{1/q}\|\psi(\cdot, s_0)\|_{L^\infty}^{1-1/q}.
\end{equation*}
\end{enumerate}
We combine all these inequalities in order to obtain the following estimation for $I_{1,B_{f_1}}$:
\begin{equation*}
I_{1,B_{f_1}}\leq C\mu f_1^{N(1-1/q)+\omega-1}\|\psi(\cdot, s_0)\|_{L^1}^{1/q}\|\psi(\cdot, s_0)\|_{L^\infty}^{1-1/q}.
\end{equation*}
\item[(ii)] \underline{Estimations for the dyadic corona $E_k$}. Let us note $I_{1,E_k}$ the integral
$$I_{1,E_k}=\int_{E_k}\|x\cdot x^{-1}(s_1)\|^{\omega-1}|v-\overline{v}_{B_{f_1}}| |\psi_+(x,s_1)|dx.$$
Since over $E_k$ we have $\|x\cdot x^{-1}(s_1)\|^{\omega-1}\leq C 2^{k(\omega-1)}f_1^{\omega-1}$ we write
\begin{eqnarray*}
I_{1,E_k}&\leq & C2^{k(\omega-1)}f_1^{\omega-1}\left(\int_{E_k}|v-\overline{v}_{B(f_12^k)}| |\psi_+(x,s_1)|dx+\int_{E_k}|\overline{v}_{B_{f_1}}-\overline{v}_{B(f_12^k)}| |\psi_+(x,s_1)|dx\right)\\
&\leq& C2^{k(\omega-1)}f_1^{\omega-1}\left(\int_{B(f_12^k)}|v-\overline{v}_{B(f_12^k)}| |\psi_+(x,s_1)|dx\right.\\
& &\qquad \qquad \qquad \qquad\left.+\int_{B(f_12^k)}|\overline{v}_{B_{f_1}}-\overline{v}_{B(f_12^k)}| |\psi_+(x,s_1)|dx\right).
\end{eqnarray*}
where $B(f_12^k)= \{x\in \mathbb{G}: \|x\cdot x^{-1}(s_1)\|\leq f_12^k\}$.\\

Now, since $v(\cdot, s_1)\in bmo(\mathbb{G})$, using the Lemma \ref{PropoBMO1} we have $|\overline{v}_{B_{f_1}}-\overline{v}_{B(f_12^k)}| \leq Ck\|v(\cdot,s_1)\|_{bmo}\leq Ck\mu$. We write 
\begin{eqnarray*}
I_{1,E_k}&\leq &C2^{k(\omega-1)}f_1^{\omega-1}\left(\int_{B(f_12^k)}|v-\overline{v}_{B(f_12^k)}| |\psi_+(x,s_1)|dx+Ck\mu\|\psi_+(\cdot,s_1)\|_{L^1}\right)\\[5mm]
&\leq & C2^{k(\omega-1)}f_1^{\omega-1}\left(\|\psi_+(\cdot,s_1)\|_{L^{a_0}} \|v-\overline{v}_{B(f_12^k)}\|_{L^{\frac{a_0}{a_0-1}}} +Ck\mu\;\|\psi_+(\cdot,s_0)\|_{L^1}\right)
\end{eqnarray*}
where we used Hölder's inequality with $1<a_0<\frac{N}{N+(\omega-1)}$ and maximum principle for the last term above. Using again the properties of $bmo$ spaces we have
$$I_{1,E_k}\leq  C2^{k(\omega-1)}f_1^{\omega-1}\left(\|\psi_+(\cdot,s_0)\|_{L^1}^{1/a_0}\|\psi_+(\cdot,s_0)\|_{L^\infty}^{1-1/a_0} |B(f_12^k)|^{1-1/a_0}\|v(\cdot,s_1)\|_{bmo} +Ck\mu\|\psi(\cdot,s_0)\|_{L^1}\right).$$
Since $\|v(\cdot,s_1)\|_{bmo}\leq \mu$ and since $1<a_0<\frac{N}{N+(\omega-1)}$, we have $N(1-1/a_0)+(\omega-1)<0$, so that, summing over each dyadic corona $E_k$, we obtain
\begin{equation*}
\sum_{k\geq 1}I_{1,E_k}\leq C \mu\left(f_1^{N(1-1/a_0)+\omega-1}\|\psi(\cdot,s_0)\|_{L^1}^{1/a_0}\|\psi(\cdot,s_0)\|_{L^\infty}^{1-1/a_0}+f_1^{\omega-1}\|\psi(\cdot,s_0)\|_{L^1}\right).
\end{equation*}
\end{enumerate}
We finally obtain the following inequalities:
\begin{eqnarray}
I_1&=& I_{1,B_{f_1}}+\sum_{k\geq 1}I_{1,E_k}\label{FinalEstimateI_1}\\
&\leq &C\mu \underbrace{f_1^{N(1-1/q)+\omega-1}\|\psi(\cdot, s_0)\|_{L^1}^{1/q}\|\psi(\cdot, s_0)\|_{L^\infty}^{1-1/q}}_{(a)}\nonumber\\
& &+ C \mu \left(\underbrace{f_1^{N(1-1/a_0)+\omega-1}\|\psi(\cdot,s_0)\|_{L^1}^{1/a_0}\|\psi(\cdot,s_0)\|_{L^\infty}^{1-1/a_0}}_{(b)} +\underbrace{f_1^{\omega-1}\|\psi(\cdot,s_0)\|_{L^1}}_{(c)}\right)\nonumber
\end{eqnarray}
Now we will prove that each of the terms $(a)$, $(b)$ and $(c)$ above is bounded by the quantity $\big(r+Ks_0\big)^{\omega-\gamma-1}$:
\begin{itemize}
\item for the first term (a) by the hypothesis on the initial data $\psi(\cdot,s_0)$ and the definition of $f_1$ given in (\ref{DefiFunctionF}) we have:
\begin{eqnarray*}
f_1^{N(1-1/q)+\omega-1}\|\psi(\cdot, s_0)\|_{L^1}^{1/q}\|\psi(\cdot, s_0)\|_{L^\infty}^{1-1/q}& \leq & \big(r+Ks_0\big)^{[N(1-1/q)+\omega-1]-\frac{\gamma}{q}-(N+\gamma)(1-1/q)}=\big(r+Ks_0\big)^{\omega-\gamma-1}.
\end{eqnarray*}
\item For the second term (b) we have, by the same arguments:
\begin{eqnarray*}
f_1^{N(1-1/a_0)+\omega-1}\|\psi(\cdot,s_0)\|_{L^1}^{1/a_0}\|\psi(\cdot,s_0)\|_{L^\infty}^{1-1/a_0} &\leq &
\big(r+Ks_0\big)^{[N(1-1/a_0)+\omega-1]-\frac{\gamma}{a_0}-(N+ \gamma)(1-1/a_0)}=\big(r+Ks_0\big)^{\omega-\gamma-1}.
\end{eqnarray*}
\item Finally, for the last term (c) we write
\begin{eqnarray*}
f_1^{\omega-1}\|\psi(\cdot,s_0)\|_{L^1}&\leq &f_1^{\omega-1}(r+Ks_0)^{-\gamma}= \big(r+Ks_0\big)^{\omega-\gamma-1}.
\end{eqnarray*}
\end{itemize}
Gathering these estimates on $(a), (b)$ and $(c)$, and getting back to (\ref{FinalEstimateI_1}) we finally obtain
$$I_1\leq C\mu \big(r+Ks_0\big)^{\omega-\gamma-1}.$$
The Lemma \ref{Lemme41} is proven.\hfill$\blacksquare$\\

\textbf{\textit{Proof of the Lemma \ref{Lemme42}.}} As for the Lemma \ref{Lemme41}, we consider $\mathbb{G}$ as the union of a ball with dyadic coronas centered on $x(s_1)$ (cf. (\ref{Decoupage})). 
\begin{enumerate}
\item[(i)] \underline{Estimations over the ball $B_{f_1}$}. We will follow closely the computations of the Lemma \ref{Lemme2}. We write:
\begin{eqnarray}
I_{2,B_{f_1}}&=&\int_{B_{f_1}}\|x\cdot x^{-1}(s_1)\|^{\omega-1}\,|\psi_+(x,s_1)|dx \leq \|\psi_+(\cdot, s_1)\|_{L^{\infty}}\int_{B_{f_1}}\|x\cdot x^{-1}(s_1)\|^{\omega-1}dx\nonumber\\
&\leq & Cf_1^{N+\omega-1}\|\psi_+(\cdot, s_0)\|_{L^\infty}.\label{Bola2}
\end{eqnarray}
\item[(ii)]  \underline{Estimations for the dyadic corona $E_k$}. Here we have
\begin{eqnarray*}
I_{2,E_k}=\int_{E_k}\|x\cdot x^{-1}(s_1)\|^{\omega-1}\, |\psi_+(x,s_1)|dx & \leq &\|\psi_+(\cdot, s_0)\|_{L^1} \underset{f_12^{k-1}<\|x\cdot x^{-1}(s_1)\|\leq f_12^k}{\sup}\|x\cdot x^{-1}(s_1)\|^{\omega-1}.\\
&\leq &C \big(2^kf_1\big)^{\omega-1} \|\psi_+(\cdot, s_0)\|_{L^1} 
\end{eqnarray*}
Since $0<\gamma<\omega<1$, summing over $k\geq 1$, we obtain
\begin{equation}\label{Coronak2}
\sum_{k\geq 1}I_{2,E_k}\leq C f_1^{\omega-1}\|\psi(\cdot, s_0)\|_{L^1}.
\end{equation}
\end{enumerate}
To finish the proof of the Lemma \ref{Lemme42} we combine (\ref{Bola2}) and (\ref{Coronak2}) and we obtain
$$I_2=I_{2,B_{f_1}}+\sum_{k\geq 1}I_{2,E_k}\leq C\left(\underbrace{f_1^{N+\omega-1}\|\psi_+(\cdot, s_0)\|_{L^\infty}}_{(d)}+\underbrace{f_1^{\omega-1}\|\psi(\cdot, s_0)\|_{L^1}}_{(e)}\right)$$
Now, we prove that the quantities $(d)$ and $(e)$ can be bounded by $\big(r+Ks_0\big)^{\omega-\gamma-1}$.
\begin{itemize}
\item For the term $(d)$ we write $f_1^{N+\omega-1}\|\psi(\cdot, s_0)\|_{L^\infty}\leq f_1^{N +\omega-1}(r+Ks_0)^{-(N+\gamma)}= \big(r+Ks_0\big)^{\omega-\gamma-1}$.
\item To treat the term $(e)$ it is enough to apply the same arguments used to prove the part $(c)$ above. 
\end{itemize}
Finally, we obtain
$$I_2=I_{2,B_{f_1}}+\sum_{k\geq 1}I_{2, E_k}\leq  C\big(r+Ks_0\big)^{\omega-\gamma-1}$$
The Lemma \ref{Lemme42} is proven.\hfill$\blacksquare$\\

\subsubsection*{2) The Height condition}
Now we write down the maximum principle for a small time $s_1$ but with a initial condition $\psi(\cdot,s_0)$, with $s_0>0$. The proof follows essentially the same ideas explained in the previous step. Indeed, since we have assumed that the Concentration condition (\ref{Concentration2}) is bounded by $(r+K(s_0+s_1))^{\omega-\gamma}$, we obtain in the same manner and with the same constants:
$$\frac{d}{ds_1}\|\psi(\cdot,s_1)\|_{L^\infty}\leq -C(r+K(s_0+s_1))^{-\frac{(\omega-\gamma)}{N+\omega}}\|\psi(\cdot,s_1)\|_{L^\infty}^{1+\frac{1}{N+\omega}}.$$
To conclude, it is enough to solve the previous differential inequality with initial data $\|\psi(\cdot,0)\|_{L^\infty}\leq (r+Ks_0)^{-(N+\gamma)}$ to obtain that $\|\psi(\cdot,s_1)\|_{L^\infty}\leq (r+K(s_0+s_1))^{-(N+\gamma)}$. 

\subsubsection*{3) The $L^1$ condition}
The $L^1$-norm condition is a direct consequence of the previous concentration condition (\ref{Concentration2}) and of the height condition (\ref{Linftyevolutionnotsmalltime}).\\

We have the estimates (\ref{Concentration2}), (\ref{Linftyevolutionnotsmalltime}) and (\ref{L1evolutionsmalltime}) and the Theorem \ref{Generalisacion} is thus proven. \hfill$\blacksquare$
\subsection{The iteration}\label{SecIterationMol}

In sections \ref{SecEvolMol1} and \ref{SecEvolMol2} we studied respectively the evolution of small molecules from time $0$ to a small time $s_0$ and from this time $s_0$ to a larger time $s_0+s_1$ and we obtained a good $L^1$ control for such molecules. It is now possible to reapply the previous Theorem \ref{Generalisacion} in order to obtain a larger time control of the $L^1$ norm. The calculus of the $n$-th iteration will be essentially the same.
\begin{Theoreme}\label{Generalisacionfin} Set $\gamma$ and $\omega$ two real numbers such that $0<\gamma<\omega<1$. Let $0< s_n\leq T$ and let $\psi(x,s_n)$ be a solution of the problem
\begin{equation}\label{Evolutionfin}
\left\lbrace
\begin{array}{rl}
\partial_{s_n} \psi(x,s_n)=& -\nabla\cdot(v\, \psi)(x,s_n)-\mathcal{J}^{1/2}\psi(x,s_n)\\[5mm]
\psi(x,0)=& \psi(x,s_{n-1})  \qquad \qquad \mbox{with } s_{n-1}>0\\[5mm]
div(v)=&0 \quad \mbox{and }\; v\in L^{\infty}([0,T];bmo(\mathbb{G}))\quad \mbox{with } \underset{s_n\in [s_{n-1},T]}{\sup}\; \|v(\cdot,s_n)\|_{bmo}\leq \mu
\end{array}
\right.
\end{equation}
If $\psi(x,s_{n-1})$ satisfies the three following conditions
\begin{eqnarray*}
\int_{\mathbb{G}}|\psi(x,s_{n-1})|\|x\cdot x^{-1}(s_{n-1})\|^\omega dx &\leq &(r+K(s_0+\cdots+s_{n-1}))^{\omega-\gamma}\\
\|\psi(\cdot, s_{n-1})\|_{L^\infty}\leq  \frac{1}{\left(r+ K(s_0+\cdots+s_{n-1})\right)^{N+\gamma }}&; &\quad \|\psi(\cdot, s_{n-1})\|_{L^1} \leq  \frac{v_n}{\big(r+K(s_0+\cdots+s_{n-1})\big)^{\gamma }}
\end{eqnarray*}
where $K=K(\mu)$ is given by (\ref{SmallConstants}) and $s_n$ is such that $(r+K(s_0+\cdots+s_n))<1$. Then for all $0< s_n\leq\epsilon r$ small, we have the following estimates
\begin{eqnarray}
\int_{\mathbb{G}}|\psi(x,s_n)|\|x\cdot x^{-1}(s_n)\|^\omega dx &\leq &(r+K(s_0+\cdots+s_n))^{\omega-\gamma}  \label{Concentration2fin}\\
\|\psi(\cdot,s_n)\|_{L^\infty}&\leq & \frac{1}{\left(r+K(s_0+\cdots+s_n)\right)^{N+\gamma }}\nonumber\\
\|\psi(\cdot,s_n)\|_{L^1} &\leq & \frac{v_N}{\big(r+K(s_0+\cdots+s_n)\big)^{\gamma }}\label{L1normestimatefin}
\end{eqnarray}
\end{Theoreme}

\begin{Remarque}
\begin{itemize}
\item[]
\item[1)] Again, since $s_n$ is small and $(r+K(s_0+\cdots+s_{n-1}))<1$, we can without loss of generality assume that $(r+K(s_0+\cdots +s_n))<1$: otherwise, by the maximum principle there is nothing to prove.
\item[2)] The new molecule's center $x(s_n)$ used in formula (\ref{Concentration2fin}) is fixed by 
\begin{equation}\label{Defpointx_sfin}
\left\lbrace
\begin{array}{rl}
x'(s_n)=& \overline{v}_{B_{f_n}}=\frac{1}{|B_{f_n}|}\displaystyle{\int_{B_{f_n}}}v(y,s_n)dy\\[5mm]
x(0)=& x(s_{n-1}).
\end{array}
\right.
\end{equation}  
And here we noted $B_{f_n}=B(x(s_n),f_n)$ with $f_n$ a real valued function given by 
\begin{equation}\label{DefiFunctionFfin}
f_n=(r+K(s_0+\cdots+s_{n-1})).
\end{equation}
Note that by remark $1)$ above we have $0<f_n<1$.
\end{itemize}
\end{Remarque}
\textit{\textbf{Proof of the Theorem \ref{Generalisacionfin}.}} The proof will follow again the same scheme: we start with the Concentration condition, we continue with the Height condition: the $L^1$ and thus the $L^q$ bound will be an easy consequence of these two estimates.
\subsubsection*{1) The Concentration condition}
Write $\Omega_n(x)=\|x \cdot x^{-1}(s_n)\|^{\omega}$ and $\psi(x)=\psi_+(x)-\psi_-(x)$, by linearity and using the positivity principle we have 
$|\psi(x,s_n)|=|\psi_+(x,s_n)-\psi_-(x,s_n)|\leq \psi_+(x,s_n)+\psi_-(x,s_n)$ and we may consider the formula:
\begin{eqnarray*}
I&=&\left|\partial_{s_n} \int_{\mathbb{G}}\Omega_n(x)\psi_+(x,s_n)dx\right|=\left|\int_{\mathbb{G}}-\nabla\Omega_n(x)\cdot x'(s_n)\psi_+(x,s_n)+\Omega_n(x)\left[-\nabla\cdot(v\, \psi_+(x,s_n))-\mathcal{J}^{1/2}\psi_+(x,s_n)\right]dx\right|
\end{eqnarray*}
Using the definition of $x'(s_n)$ given in (\ref{Defpointx_sfin}) and replacing $\Omega_n(x)$ by $\|x \cdot x^{-1}(s_n)\|^{\omega}$ in the first integral we obtain
\begin{equation}\label{Estrellafin}
I\leq c \underbrace{\int_{\mathbb{G}}\|x \cdot x^{-1}(s_n)\|^{\omega-1}|v-\overline{v}_{B_f}| |\psi_+(x,s_n)|dx}_{I_1} + c\underbrace{\int_{\mathbb{G}}|\mathcal{J}^{1/2}\Omega_n(x)||\psi_+(x,s_n)|dx}_{I_2}.
\end{equation}
We will study each of the integrals $I_1$ and $I_2$ in the next lemmas:
\begin{Lemme}\label{Lemme41fin} For integral $I_1$ we have $I_1\leq C \mu\big(r+K(s_0+\cdots+s_{n-1})\big)^{\omega-\gamma-1}$.
\end{Lemme}
\textit{\textbf{Proof}.} It is enough to repeat the same steps of the previous Lemma \ref{Lemme41}, just consider $\mathbb{G}=B_{f_n}\cup\bigcup_{k\geq 1}E_k$ where
\begin{eqnarray}\label{Decoupagefin}
B_{f_n}= \{x\in \mathbb{G}: \|x\cdot x^{-1}(s_n)\|\leq f_n\}, \qquad E_k= \{x\in \mathbb{G}: f_n2^{k-1}<\|x\cdot x^{-1}(s_n)\|\leq f_n 2^{k}\} \quad \mbox{for } k\geq 1.
\end{eqnarray}
In order to obtain the desired inequality, use exactly the same arguments, the maximum principle and the hypothesis of Theorem  \ref{Generalisacionfin}.\hfill$\blacksquare$
\begin{Lemme}\label{Lemme42fin}
For integral $I_2$ in inequality (\ref{Estrellafin}) we have the following estimate
\begin{equation*}
I_2=\int_{\mathbb{G}}|\mathcal{J}^{1/2}\Omega_n(x)||\psi_+(x,s_n)|dx\leq C\big(r+K(s_0+\cdots+s_{n-1})\big)^{\omega-\gamma-1}.
\end{equation*}
\end{Lemme}
\textbf{\textit{Proof.}} As for Lemma \ref{Lemme41fin}, we consider $\mathbb{G}$ as the union of a ball with dyadic coronas centered on $x(s_n)$ (cf. (\ref{Decoupagefin})). It is then enough to repeat the corresponding estimates of the $s_1$-case given in Lemma \ref{Lemme42}. \hfill$\blacksquare$\\

Now using the Lemmas \ref{Lemme41fin} and \ref{Lemme42fin} and getting back to the estimate (\ref{Estrellafin}) we have
\begin{equation}\label{FinalEstimatefin}
\left|\partial_{s_n} \int_{\mathbb{G}}\Omega_n(x)\psi_+(x,s_n)dx\right| \leq  C \left(\mu+1\right)\big(r+K(s_0+\cdots+s_{n-1})\big)^{\omega-\gamma-1}
\end{equation}
This estimation is compatible with the estimate (\ref{Concentration2fin}) for $0\leq s_n\leq \epsilon r$ small enough. Indeed, we can write
$\phi=(r+K(s_0+\cdots+s_n))^{\omega-\gamma}$ and we linearize this expression with respect to $s_n$:
$$\phi \thickapprox(r+K(s_0+\cdots+s_{n-1}))^{\omega-\gamma}\left(1+K(\omega-\gamma)\frac{s_n}{(r+K(s_0+\cdots+s_{n-1}))}\right)$$
Taking the derivative of $\phi$ with respect to $s_n$ we have $\phi' \thickapprox K(\omega-\gamma) \big(r+K(s_0+\cdots+s_{n-1})\big)^{\omega-\gamma-1}$ and with the condition (\ref{SmallConstants}) on $K(\omega-\gamma)$ we obtain that (\ref{FinalEstimatefin}) is bounded by $\phi'$ and we have proven the Concentration condition.
\subsubsection*{2) The Height condition}
Since we have that Concentration condition (\ref{Concentration2fin}) is bounded by $(r+K(s_0+\cdots+s_n))^{\omega-\gamma}$, following the previous computations we obtain in the same manner and with the same constants:
$$\frac{d}{ds_n}\|\psi(\cdot,s_n)\|_{L^\infty}\leq -C(r+K(s_0+\cdots+s_n))^{-\frac{(\omega-\gamma)}{N+\omega}}\|\psi(\cdot,s_n)\|_{L^\infty}^{1+\frac{1}{N+\omega}}.$$
Solving this differencial inequality we obtain $\|\psi(\cdot,s_n)\|_{L^\infty}\leq (r+K(s_0+\cdots+s_n))^{-(N+\gamma)}$.
\subsubsection*{3) The $L^1$-norm estimate}
Again this is a direct consequence of the Concentration condition and of the previous Height condition.\\[5mm]
Theorem \ref{Generalisacionfin} is completely proven. \hfill$\blacksquare$
\subsubsection*{End of the proof of Theorem \ref{TheoL1control}} 
We see with the Theorem \ref{SmallGeneralisacion} that is possible to control the $L^1$ behavior of the molecules $\psi$ from $0$ to a small time $s_0$, from time $s_0$ to time $s_1$ with Theorem \ref{Generalisacion}, and by iteration from time $s_{n-1}$ to time $s_n$ with Theorem \ref{Generalisacionfin}. We recall that we have $s_i\sim \epsilon r$ for all $0\leq i\leq n$, so the bound obtained in (\ref{L1normestimatefin}) depends mainly on the size of the molecule $r$ and the number of iterations $n$. \\

We observe now that the smallness of $r$ and of the times $s_0,...,s_n$ can be compensated by the number of iterations $n$ in the following sense: fix a small $0<r<1$ and iterate as explained before. Since each small time $s_0,...,s_n $ is of order $\epsilon r$, we have $s_0+\cdots+s_{n}\sim n \epsilon r$. Thus, we will stop the iterations as soon as $ n r\geq T_0$. 

Of course, the number of iterations $n=n(r)$ will depend on the smallness of the molecule's size $r$, and more specifically it is enough to set $n(r)\sim \frac{T_0}{r}$ in order to obtain this lower bound for $nr$.

Proceeding this way we will obtain $\|\psi(\cdot,s_n)\|_{L^1}\leq C T_0^{-\gamma}<+\infty$, for all molecules of size $r$. Note in particular that, once this estimate is available, for bigger times it is enough to apply the maximum principle.\\

Finally, and for all $r>0$, we obtain after a time $T_0$ a $L^1$ control for small molecules and we finish the proof of Theorem \ref{TheoL1control} since the $L^q$ control can be easily deduced from the $L^1$ and $L^\infty$ bounds. \hfill$\blacksquare$

\quad\\[5mm]

\begin{flushright}
\begin{minipage}[r]{80mm}
Diego \textsc{Chamorro}\\[3mm]
Laboratoire d'Analyse et de Probabilités\\ 
Université d'Evry Val d'Essonne\\[2mm]
23 Boulevard de France\\
91037 Evry Cedex\\[2mm]
diego.chamorro@univ-evry.fr
\end{minipage}
\end{flushright}

\end{document}